%% file: main.tex
\documentclass{article}
\usepackage[utf8]{inputenc}
\usepackage{authblk}

\title{The Maintenance Scheduling and Location Choice Problem for Railway Rolling Stock}
\author{Jordi Zomer$^1$, Nikola Bešinović$^1$, Mathijs M. de Weerdt$^2$, Rob M.P. Goverde$^1$}
\affil{Delft University of Technology, Delft, The Netherlands \\
\small $^1$ 
Faculty of Civil Engineering and Geosciences, Department of Transport \& Planning
\\ $^2$ Faculty of Electrical Engineering, Mathematics and Computer Science, Department of Software and Computer Technology.

}

\date{March 2021}

\usepackage{subfiles}

\usepackage[utf8]{inputenc}

\usepackage{heuristica}
\usepackage[heuristica]{newtxmath}
\usepackage{fullpage}

\usepackage{longtable}
\usepackage{booktabs, ctable}

\usepackage{float}
\usepackage{placeins}
\usepackage{array}
\newcolumntype{P}[1]{>{\raggedright\arraybackslash}p{#1}}
\newcolumntype{C}[1]{>{\centering\arraybackslash}p{#1}}
\usepackage[font={small,it}]{caption}

\usepackage{lipsum}

\usepackage{enumitem}
\usepackage{graphicx}
\usepackage{subcaption}
\usepackage{lscape}

\usepackage{amsmath}

\usepackage[official]{eurosym}

\usepackage{pdfpages}

\usepackage{algorithm}
\usepackage{algpseudocode}

\usepackage{hyperref}

\usepackage{natbib}
\usepackage{url}
\begin{document}

\maketitle
\input{00_Abstract}
\input{01_Introduction} 
\input{02_Literature_review}
\input{03_MLCP}
\input{04_APP}
\input{05_MSLCP}
\input{06_Results}
\input{07_Conclusion}

\newpage \bibliography{References}

\end{document}

%% file: 00_Abstract.tex


\begin{abstract}
Due to increasing railway use, the capacity at railway yards and maintenance locations is becoming limiting. Therefore, the scheduling of rolling stock maintenance and the choice regarding optimal locations to perform maintenance is increasingly complicated.
This research introduces a Maintenance Scheduling and Location Choice Problem (\textbf{MSLCP}). It simultaneously determines maintenance locations and maintenance schedules of rolling stock, while it also considers the available capacity of maintenance locations, measured in the number of available teams.
To solve the \textbf{MSLCP}, an optimization framework based on Logic-Based Benders' Decomposition (LBBD) is proposed by combining two models, the Maintenance Location Choice Problem (\textbf{MLCP}) and the Activity Planning Problem (\textbf{APP}), to assess the capacity of a \textbf{MLCP} solution. Within the LBBD, four cut generation procedures are introduced to improve the computational performance: a naive procedure, two heuristic procedures and the so-called min-cut procedure that aims to exploit the specific characteristics of the problem at hand. The framework is demonstrated on a realistic scenarios from the Dutch railways. It is shown that the best choice for cut generation procedure depends on the objective: when aiming to find a good but not necessarily optimal solution, the min-cut procedure performs best, whereas when aiming for the optimal solution, one of the heuristic procedures is the preferred option. The techniques used in the current research are new to the current field and offer interesting next research opportunities.

\end{abstract}


%% file: 01_Introduction.tex


\section{Introduction} \label{sec:introduction}


In many countries, rail transport is increasingly important. This is for example visible in The Netherlands, where the total number of passenger kilometers increased with more than 30\% since 2018 \citep{UIC2018paxkm}.
In order for a railway network to function properly, the rolling stock 
(i.e. locomotives, passenger wagons and freight wagons, multiple units) 
that operates on the railway network needs to receive maintenance on a regular basis. 
The aim of maintenance is to ensure that the rolling stock that operates on the network remains available (to ensure a reliable train service), safe and comfortable for passengers \citep{dinmohammadi2016risk}. To this end, 
maintenance activities can be divided into two categories: regular maintenance, corresponding to the maintenance activities with higher frequencies (every 1 to 14 days) and shorter duration (1-3 hours) which can be performed whenever the unit has a planned standstill, and heavy maintenance, corresponding to maintenance types with lower frequencies (every several months or less) and longer duration (up to several days) during which the unit is taken out of service (see e.g.\ \citet{Andres2015MaintenanceNetworks}).  The current work focuses on regular maintenance in particular.
In general, the maximum interval between consecutive maintenance activities is governed by strict rules that are imposed by railway authorities. 

Maintenance activities are carried out at so-called maintenance locations, which are railway yards with maintenance facilities, spread over the network. The number of personnel stationed at a location is the operator's decision and an important determinant of the capacity of a maintenance location. In particular, a distinction can be made between daytime operations (i.e. a location is opened during daytime) and nighttime operations (i.e. a location is opened during nighttime). For example, in The Netherlands, maintenance is usually carried out during nighttime. 

The increasing use of the capacity of the railway network leads to two issues for rolling stock maintenance. First, the complexity of this scheduling process -- which is traditionally performed manually -- is increasing. This raises the need for tools that automate the maintenance scheduling process. Second, the use of the capacity of maintenance locations during nighttime is under pressure and reaching its capacity. As a result, 
a railway operator can consider to perform more maintenance activities during daytime, as is the case in The Netherlands for example.
These issues were originally addressed by \citet{zomer2020MLCP}, providing a model for the \textit{Maintenance Location Choice Problem} (\textbf{MLCP}). However, they assumed  an unlimited capacity of maintenance locations and do not provide a maintenance schedule. 
The capacity of maintenance locations
is a challenging factor to incorporate,
as it typically depends on the optimal planning of all maintenance activities (which is not readily available).  

This research introduces a new mathematical problem formulation, the Maintenance Scheduling and Location Choice Problem (\textbf{MSLCP}), which extends the \textbf{MLCP}. For a given rolling stock circulation, it determines an optimal maintenance schedule and an optimal maintenance location choice while including capacity constraints of maintenance locations. In order to do so, it introduces a separate problem to assess the capacity of a \textbf{MLCP} solution, called the \textit{Activity Planning Problem} (\textbf{APP}). 
To solve the MSLCP, an optimization framework based on Logic-Based Benders' Decomposition (LBBD) is proposed by combining two models, the \textbf{MLCP} and \textbf{APP}. 
Within the LBBD, four cut generation procedures are introduced: a naive procedure, two heuristic procedures, and lastly the so-called min-cut procedure which uses the specific structure of the problem at hand. The performance of \textbf{MSLCP} is demonstrated on a realistic case from the Dutch railways. 

The contribution of this research is threefold.
First, it extends the model proposed by \citet{zomer2020MLCP}, by introducing maintenance location capacity constraints and thereby making the model capable of delivering a complete maintenance schedule. Second, it provides an efficient method to assess the required capacity of a maintenance schedule which has as an additional benefit that it can be used to quickly provide rolling stock dispatchers with a maintenance activity planning during operations. Third, it proposes an advanced solution strategy for the inclusion of the capacity constraints based on Logic-Based Benders' Decomposition. 

The remainder of this paper is structured as follows. Section~\ref{sec:litrev} summarizes the most important existing literature. 
Section~\ref{sec:MLCP} restates the \textbf{MLCP} based on \citet{zomer2020MLCP}. Section~\ref{sec:app} formulates the \textbf{APP}, allowing to assess the capacity of any \textbf{MLCP} solution. Section~\ref{sec:MSLCP} formulates the \textbf{MSLCP}, in which the \textbf{MLCP} and \textbf{APP} are the main building blocks, and provides a solution strategy. Section~\ref{sec:results} provides results for the \textbf{MSLCP} and Section~\ref{sec:conclusion} gives the main conclusions.


%% file: 02_Literature_review.tex

\section{Literature review} \label{sec:litrev}
The current work considers rolling stock maintenance scheduling as well as rolling stock maintenance location choice. This section aims to identify the contributions of the current work to the literature and to obtain insights in the methodologies and techniques used in related research.

Section~\ref{sec:litrev_sched} discusses relevant scientific literature on rolling stock maintenance scheduling. Section~\ref{sec:litreview_maintlocation} discusses papers on rolling stock maintenance location choice. 
In addition, these two sections include some corresponding literature from the field of aviation, which is relevant due to the systematic similarities with rolling stock maintenance scheduling. Finally, Section \ref{sec:litrev_contribution} summarises reviewed papers and states existing gaps.

\subsection{Maintenance scheduling}~\label{sec:litrev_sched}
\citet{Herr2017F.Scheduling} considered a problem in which rolling stock units need to be assigned to train trips such that maintenance constraints are satisfied. They proposed a MIP model and the objective that they used is to schedule maintenance as late as possible, thereby making optimal use of the total allowable interval between maintenance activities. 

Just as \citet{Herr2017F.Scheduling}, \citet{Andres2015MaintenanceNetworks} considered the problem of assigning rolling stock units to train trips. They used an aggregated space-time network in which the nodes are trip arrival times or trip departure times with the corresponding location. A MIP model that minimizes total operating costs was designed and a column generation approach was used to solve the problem in reasonable time.

\citet{Maroti2007MaintenanceModel} considered a problem regarding heavy maintenance. They proposed a model to make modifications to the regular rolling stock circulation to route rolling stock units to maintenance locations and formulated it as an integer programming problem. In situations where only one rolling stock unit needs to be rerouted, this formulation provides the optimal solution; in situations where multiple rolling stock units need to be rerouted the formulation is used within a heuristic framework.

\citet{Wagenaar2015MaintenanceRailways} proposed with a model that reschedules rolling stock circulation after disruptions taking into account the current maintenance planning. They based their models on the composition model, which assigns rolling stock units to train trips. They came up with three models that have comparable performance, dependent on the problem size.

Related problems were addressed in the area of aviation, for example by \citet{clarke1997aircraft} and \citet{gopalan1998aircraft}, who aimed to assign specific aircraft to each flight from a given set of flights, and \citet{Sarac2006ARouting}, who developed a model that solves the aircraft maintenance scheduling problem including maintenance constraints in an operational context. 



\subsection{Maintenance location choice} \label{sec:litreview_maintlocation}
\citet{Tonissen2019MaintenanceUncertainty} aimed at locating the maintenance facilities in the railway network. They came up with models that determine optimal maintenance locations under line and fleet planning that is subject to uncertainty or change. They proposed two-stage stochastic mixed integer programming models, in which the first stage is to open a facility, and in the second stage to minimize the routing cost for the first-stage location decision for each line plan scenario.

\citet{Tonissen2018EconomiesStock} built on \citet{Tonissen2019MaintenanceUncertainty} by including recovery costs of maintenance location decisions, unplanned maintenance, multiple facility sizes and economies of scale (providing that a location twice as big is not twice as expensive). Since, as a result, the second-stage problem becomes NP-hard, an algorithm was provided with the aim to avoid having to solve the second stage for every scenario. 

\citet{Canca2018TheSystems} considered the simultaneous rolling stock allocation to lines and choice for depot locations in a rail-rapid transit context. They proposed a MIP formulation which appeared hard to solve. Therefore they proposed a three-step heuristic approach determining first the minimum number of vehicles needed for each line, subsequently the routes of rolling stock on each line, and lastly the circulation of rolling stock on lines over multiple days together with the depot choice.

\citet{zomer2020MLCP} introduced the Maintenance Location Choice Problem (\textbf{MLCP}). To solve it, the authors developed a Mixed Integer Linear Programming model taking a rolling stock circulation as input, and provided for this rolling stock circulation an optimal maintenance location choice that minimize the total number of maintenance activities during nighttime, thereby reducing the pressure on maintenance locations during nighttime. 
However, they did not include the capacity of maintenance locations, nor determine exact maintenance schedules, i.e.   maintenance activities are assigned to maintenance opportunities, which are longer time windows in which maintenance has to take place at some moment, and do not consider actual moment when maintenance has to be performed.



Some related research can be found in the area of aviation. Examples are the works by \citet{Feo1989FlightPlanning} and \citet{Gopalan2014TheProblem}, who consider the problem of assigning aircraft to flights and simultaneously determining maintenance locations and introduce various heuristics to solve the problem.

\subsection{Current work}~\label{sec:litrev_contribution}
In Table~\ref{tab:lit_overview}, the discussed literature is classified in several categories. It shows for each paper whether it was written in the aviation (A) or in the railway (R) context, whether it considered the allocation of mobile units (MU, i.e. rolling stock units or aircraft) to trips, whether it considered maintenance constraints, whether it created an explicit maintenance schedule for every (relevant) MU and whether it considered facility location choice optimization. 

\begin{table}[H]
    \centering
    \small
\begin{tabular}{lC{2cm}C{2cm}C{2cm}}
    \FL
       & \multicolumn{1}{p{2cm}}{MU \newline allocation}  &  \multicolumn{1}{p{2cm}}{Maintenance scheduling}  & \multicolumn{1}{p{2cm}}{ Location choice} 
     \ML
     \citet{Herr2017F.Scheduling} & x & x &  \\
     \citet{Andres2015MaintenanceNetworks} & x & x & \\ 
     \citet{Maroti2007MaintenanceModel} & x & x & \\ 
     \citet{Wagenaar2015MaintenanceRailways} & x & x & \\ 
     \citet{Tonissen2019MaintenanceUncertainty} & & & x \\
     \citet{Tonissen2018EconomiesStock} & & & \\
     \citet{Canca2018TheSystems} & x & & x \\
     \citet{zomer2020MLCP} & &  & x \\
     \textit{Current} & & x & x
     \LL
\end{tabular}
    \caption{Overview of the literature discussed in Section~\ref{sec:litrev}.}
    \label{tab:lit_overview}
\end{table}

This literature review indicates that several aspects have not been addressed in the currently existing literature. 
First, although variants of problems relating to rolling stock maintenance location and maintenance scheduling have been investigated independently, this joint problem has not been tackled. That is, no research has aimed to determine optimal opening of maintenance locations and simultaneously find an optimal maintenance schedule for a given rolling stock circulation. 
Second, although some papers do consider some type of a constraint for the available capacity at maintenance locations, such constraints are typically rather general and ignore many practical aspects. This research instead accurately determines the capacity of maintenance locations 
by first scheduling maintenance activities optimally. Moreover, the capacity is measured as the minimal number of maintenance teams necessary to fulfil a certain schedule, which is also new.
Third, the current paper delivers an exact maintenance schedule, which provides operators at maintenance location with exact moments when each rolling stock unit needs to be maintained and by which maintenance team.



%% file: 03_MLCP.tex


\section{Maintenance Location Choice Problem (MLCP)} \label{sec:MLCP}
This section summarizes the mathematical model of  the \textit{Maintenance Location Choice Problem} (\textbf{MLCP}) \cite{zomer2020MLCP}. For more detailed explanations of the model and the computational experiments, the reader may resort to \citet{zomer2020MLCP}. 

The following notation is used for the parameters of the model. Let $I$ be the set of rolling stock units, $T \in \mathbb{R}$ the planning horizon in hours and $L$ the set of potential maintenance locations. 
The rolling stock circulation is assumed to be given. 
A \textit{maintenance opportunity} (MO) occurs when a rolling stock unit is standing still at a potential maintenance location. 
Let $J_i \equiv \{1,..., \overline{J_i}\}$ denote the MOs for rolling stock unit $i \in I$. The location of a rolling stock unit $i$ at MO $j \in J_i$ is denoted by $l_{ij} \in L$. The start time of MO $j \in J_i$ is denoted by $s_{ij} \in \mathbb{R}$ and the end time by $e_{ij} \in \mathbb{R}$.

Let $d_{ij}=1$ indicate that an MO occurs during daytime:
\begin{align}    
d_{ij} = \begin{cases}
1 & \text{ if } \delta^D \leq e_{ij} \text{ mod } 24 < \delta^N \notag \\
0 & \text{ else} \end{cases},
\end{align}
where $\delta^D$ is the time daytime maintenance starts and $\delta^N$ the time nighttime maintenance starts. Unless stated otherwise, $\delta^D = 7.00$ and $\delta^N = 19.00$.
Let $K$ be the set of \textit{maintenance types}, $K \equiv \{1, ..., \overline{K}\}$. For each maintenance type $k \in K$, let $v_k \in \mathbb{R^+}$ be its duration in hours and let $o_k \in \mathbb{R^+}$ be the maximum interval between two consecutive maintenance activities of maintenance type $k$ in hours.

The decision variable $y_l^D \in \{0,1\}$ to \emph{open a potential maintenance location during daytime} is 1 if location $l \in L$ is available for daytime maintenance. Similarly $y_l^N \in \{0,1\}$ is equal to 1 if location $l \in L$ is available for nighttime maintenance. The number of potential maintenance locations that can be opened during daytime is restricted by $L_{max}^D$:
\begin{align}
    \sum_{l \in L} y_l^D & \leq L_{max}^D. \label{eq:M11-c5}
\end{align}

The \emph{assignment decisions} of maintenance activities to maintenance opportunities is encoded by $x_{ijk} \in \{0,1\}$, which is 1 if maintenance of type $k$ is performed to rolling stock unit $i \in I$ at MO $j \in J_i$, and 0 otherwise. It is required that the total time available at MO $j$ is not exceeded: 
\begin{align}
\sum_{k \in K}x_{ijk}v_k &\leq e_{ij} - s_{ij} & \forall i \in I, j \in J_i . \label{eq:M11-c4}
\end{align}
Furthermore, an MO $j$ can only be used if the corresponding location is open at the moment of the MO. Therefore, if $d_{ij} = 0$ and $y_{l_{ij}}^N = 0$ then $\forall_{k \in K} x_{ijk} = 0$, and similarly if $d_{ij} = 1$ and $y_{l_{ij}}^D = 0$ then $\forall_{k \in K} x_{ijk} = 0$, which is encoded in a single linear constraint as follows:
\begin{align}
x_{ijk} &\leq y^D_{l_{ij}} \cdot d_{ij} + y^N_{l_{ij}} \cdot (1-d_{ij}) & \hspace{0.5cm}  \forall i \in I, j \in J_{i}, k \in K. \label{eq:M11-c3}
\end{align}

Finally, the intervals between two successive maintenance activities $j$ and $j'$ of the same type $k$ should be at most $o_k$ apart.
This is modeled as follows: if $x_{ijk}=1$ (and ${e_{ij}+o_k \leq T}$) then $\exists j' \in V_{ijk} : x_{ij'k} = 1$, where $V_{ijk} = \{ p \in J_i : e_{ij} < s_{ip} \leq e_{ij} + o_k \}$.
For a correct start, let $b_{ik}$ be the time since the last maintenance activity of type $k$ for rolling stock unit $i$ at the start of the planning horizon, let $V_{i0k} = \{ p \in J_i : s_{ip} \leq o_k + b_{ik} \}$ and
\begin{align}
    1 &\leq \sum_{p \in V_{i0k}} x_{ipk} & \forall i \in I, k \in K \label{eq:M11-c1} \\
    x_{ijk} &\leq \sum_{p \in V_{ijk}} x_{ipk} & \forall i \in I, j \in J_i, k \in K : e_{ij} + o_k \leq T. \label{eq:M11-c2}
\end{align}

The model aims to find $x_{ijk}$ and $y_l$ satisfying the constraints (\ref{eq:M11-c5}) to (\ref{eq:M11-c2}) that minimize number of maintenance activities during the night:
\begin{align}
    \min \sum_{i \in I} \sum_{j \in J_i} \sum_{k \in K} x_{ijk} (1-d_{ij}) + \varepsilon \sum_{i \in I} \sum_{j \in J_i} \sum_{k \in K} x_{ijk}. \label{eq:M11-obj}
\end{align}
The second term penalizes every maintenance activity with an arbitrarily small penalty cost $\varepsilon$ in order to avoid unnecessary maintenance activities being performed.

%% file: 04_APP.tex

\section{Activity Planning Problem (APP)} \label{sec:app}
The \textbf{MLCP} delivers an assignment of maintenance activities to maintenance opportunities. Maintenance activities are not explicitly scheduled accurate to the minute, but may be performed anytime within the MO. This alone does not allow to directly determine the required number of maintenance teams to effectuate the maintenance schedule.  

Therefore, the \textit{Activity Planning Problem} (\textbf{APP}) is defined. The input of the \textbf{APP} is a set of jobs that need to be performed and a maximum number of maintenance teams available. A job represents the activities that need to be performed on one rolling stock unit during a specified maintenance opportunity. A job may contain one maintenance activity of a specific maintenance type, but can also contain multiple maintenance activities of different maintenance types. 

An important reason why different maintenance activities on the same rolling stock unit are grouped into jobs is that, in practice, maintenance activities of different types on the same rolling stock unit often cannot be performed simultaneously, e.g. external cleaning and wheels inspection. To ensure this, the slightly stricter assumption is made that maintenance activities on one rolling stock unit need to be performed subsequently and uninterruptedly. This assumption is deemed acceptable in practice. It also 
simplifies the model as it is not necessary to include separate, complicating constraints to 
prohibit that maintenance activities of different types on the same rolling stock units are performed simultaneously. 

The most important output of the \textbf{APP} is the number of teams necessary to perform the given set of jobs, which is a measure of the required capacity of a \textbf{MLCP} solution (or no output if this number exceeds the maximum number of maintenance teams available). This measure is essential for the development of the \textbf{MSLCP} in Section~\ref{sec:MSLCP}. Additionally, the \textbf{APP} gives the corresponding optimal activity planning, defining the start and end times of each job, which has useful practical applications as well. 

The \textbf{APP} can be applied on each individual maintenance shift to determine the required number of maintenance teams. A maintenance shift is a period of time for which a planning is made. The current research assumes each day contains two maintenance shifts: a daytime shift and a nighttime shift. 
Therefore, the output of APP are the shift plan and the required capacity.

The \textbf{APP} shows similarities with the class of \textit{Parallel Machine Scheduling Problems}, as addressed by for example~\citet{kravchenko2009minimizing}.

Section \ref{sec:activityplanning_APPchar} gives mathematical notation of APP.
Section \ref{sec:activityplanning_MSLCPtoAPP} explains the interaction between \textbf{MLCP} and \textbf{APP} including assignment of maintenance activities to shifts, setting release and deadline times and job durations.
Section \ref{sec:APPmodel} gives the model formulation.

\subsection{\textbf{APP} mathematical notation} \label{sec:activityplanning_APPchar}

\paragraph{Jobs}Let $Q$ be a given set of jobs that need to be scheduled, and for each job $q \in Q$ let the release time $r_q \in \mathbb{R}$, the deadline time $t_q \in \mathbb{R}$ and the duration $v_q \in \mathbb{R}$ be given. For each shift at a location, one job for every maintenance opportunity that has maintenance activities assigned at this location during this shift. The duration of the job is the sum of the duration of these main activities. 
The release time and deadline are based on the start and end time of the respective maintenance opportunity (more details in Section \ref{sec:activityplanning_MSLCPtoAPP}).
When constructing APP instances it is ensured that for each job, the time between the release and deadline of the job is larger than or equal to the job duration. 

\paragraph{Teams}
Each job needs to be performed by one and only one maintenance team. The team works on this job uninterruptedly, i.e. the job cannot be split into multiple separate parts (meaning preemption is not allowed). Let $\overline{N}$ be the maximum number of available maintenance teams and define $N = \{1, ..., \overline{N} \}$ to be the set of maintenance teams.

\paragraph{Scheduling}
The maintenance jobs are assigned to maintenance teams, and the start time of each maintenance job is determined. The end time of the job is then automatically determined by adding the job duration to the start of the maintenance job. The start time should be such that it is after the release time of a job, and such that the end time is before the deadline of a job. 

The current formulation of the \textbf{APP} uses so-called \textit{moments}. A moment represents the opportunity of a maintenance team to start a job. This is a construct used to model the \textbf{APP} as a linear problem. Each team has a set of moments available, corresponding to the maximum number of jobs that they can perform. To any moment, a job can be assigned. If a job is assigned to a moment, the start time of this particular moment is associated to the start time of the corresponding maintenance job. The introduction of the concept of moments allows to model a sequential planning, by requiring that if a job is assigned to moment $m$, moment $m+1$ can start only after the job assigned to moment $m$ is finished. 

Let $\overline{M}$ be the number of moments available per team and define $M = \{1, ..., \overline{M} \}$ to be the set of moments. Note that the maximum number of moments used by a team occurs when a team is continually occupied with maintenance activities of the shortest duration for the entire length of the maintenance shift. A sufficiently large $\overline{M}$ can thus be obtained by dividing the total time available in a maintenance shift over the minimum time required for each maintenance job. Moreover, the number of moments necessary never exceeds the total number of jobs. Based on these two indications, Equation~(\ref{eq:capacity_optimalM}) gives an appropriate value for $\overline{M}$ that is used in the current research. 
\begin{align}
     \overline{M} = \min \left( \left\lceil \frac{\delta^N - \delta^D}{\min_{k \in K} v_k} \right\rceil, |Q| \right) \label{eq:capacity_optimalM}
\end{align}

\paragraph{Objective}
The objective is to minimize the number of available maintenance teams necessary.

\subsection{From \textbf{MLCP} output to \textbf{APP} input} \label{sec:activityplanning_MSLCPtoAPP}
\subsubsection{Assignment of maintenance activities to maintenance shifts}
The current research assumes maintenance shifts in two time windows: the daytime maintenance shift between 7.00 and 19.00 and the nighttime maintenance shift between 19.00 and 7.00. 

Unique maintenance shifts are characterised by a combination of maintenance location, time window (i.e. daytime or nighttime) and reference day (i.e. the day when the maintenance shift starts).  An example of a unique maintenance shift would be \textit{the night shift in Amsterdam on day 3}, meaning the shift that starts in Amsterdam at 19.00 on day 3 and ends in Amsterdam at 07.00 on day 4.

The following procedure is used to determine to what maintenance shift an MO belongs.
\begin{itemize}
    \item Suppose an MO is classified as a daytime MO. Then, by the definition of daytime MOs, it is clear that the entire MO is contained within the daytime time window. The reference day is therefore equal to the end time of the MO and it belongs to the daytime maintenance shift of that particular day. 
    \item Suppose an MO is classified as a nighttime MO. Note that this does not necessarily mean that the end time is during nighttime (for example, an MO starting during nighttime and ending during daytime is still classified as a nighttime MO). For nighttime MOs, it seems reasonable to assign these MOs to the last nighttime maintenance shift that it was in. In other words, if the end time is between 0.00 and 19.00, it is classified as an MO during the nighttime shift with a reference day at the \textit{previous day}; if, on the other hand, the end time is between 19.00 and 0.00, this last maintenance shift is the nighttime maintenance shift with reference day on the \textit{current day}.
\end{itemize}

\begin{figure}[H]
    \centering
    \includegraphics{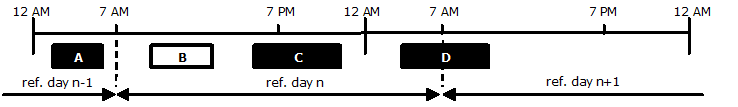}
    \caption{Assignment of MOs to shifts. This figure presents four MOs (A, B, C and D), of which B is classified as daytime MO and A, C and D are classified as nighttime MOs in the \textbf{MLCP}.}
    \label{fig:capacity_mostoshifts}
\end{figure}

An example is found in Figure~\ref{fig:capacity_mostoshifts}. This figure presents four MOs with their start and end time. Based on the above described procedure, these MOs can be assigned to maintenance shifts: MO $A$ is assigned to the nighttime shift of reference day $n-1$, MO $B$ is assigned to the daytime shift of reference day $n$ and MOs $C$ and $D$ are assigned to the nighttime shift of reference day $n$.

\subsubsection{Release and deadline times}
For each of the maintenance jobs considered by the \textbf{APP}, a release time and a deadline time need to be specified. 
It must be noted that in most cases, MOs are contained in either the daytime shift or the nighttime shift. In these cases, the release time is equal to the start of the MO and the deadline time is equal to the end of the MO. However, there are also MOs that are not fully contained in the corresponding maintenance shift, such as MOs $C$ and $D$ in Figure~\ref{fig:capacity_mostoshifts}. Still, they are assigned to nighttime maintenance shift $n$ and therefore need to be performed in this shift. 

In order to make sure that maintenance activities are performed as much as possible in the maintenance shift that they were assigned to, the following rules are used to determine the release times.
\begin{itemize}
    \item If a maintenance activity takes place in a daytime MO, then its release time is equal to the start of the corresponding MO.
    \item If a maintenance activity takes place in a nighttime MO and the start of the MO is after the start of the maintenance shift, then the release time of the maintenance job is equal to the start of the corresponding MO.
    \item If a maintenance activity takes place in a nighttime MO and the start of the MO is before the start of the maintenance shift, then the release time of the maintenance job is set to the start of the maintenance shift (usually 19.00). There is one exception to this rule: when, by setting the release time to 19.00, the time available for maintenance (i.e. between the end of the MO and 19.00) is less than the duration of the maintenance, then the release time is set to end time minus the total duration of maintenance in this job. 
\end{itemize}

\noindent A similar, symmetric set of rules prevails for the determination of the deadline moment.
\begin{itemize}
    \item If a maintenance activity takes place in a daytime MO, then its deadline time is equal to the end of the MO.
    \item If a maintenance activity takes place in a nighttime MO and the end of the MO is before the end of the maintenance shift, then the release time of the maintenance job is equal to the start of the MO.
    \item If a maintenance activity takes place in a nighttime MO and the end of the MO is after the end of the maintenance shift, then the deadline time of the  maintenance job is set to the end of the maintenance shift (usually 07.00). There is one exception to this rule: when, by setting the deadline time to 07.00, the time available for maintenance (i.e. between the start of the MO and 07.00) is less than the duration of the maintenance, then the deadline time is set to start time plus the total duration of maintenance in this job.
\end{itemize}

\subsection{\textbf{APP} model formulation} \label{sec:APPmodel}
With the notation and definitions above, the APP can now be defined.
Let $z_{nmq} \in \{0,1\}$ be a binary variable that signifies whether moment $m$ for team $n$ is associated to job $q$, where $z_{nmq} = 1$ if and only if team $n$ at moment $m$ processes job $q$. Let $s_{nm} \in \mathbb{R}$ be the start time of the moment $m$ for team $n$, where $y_n \in \{0, 1 \}$ be a binary variable that signifies whether team $n$ is active or not: let $y_n = 1$ if team $n$ is used for this schedule.

The \textbf{APP} model is then formulated as follows. 

\begin{align}
    \min \sum_{n \in N} y_n \label{eq:APP-obj}
\end{align}
subject to
\begin{align}
    \sum_{q \in Q} z_{nmq} r_q & \leq s_{nm} \leq \sum_{q \in Q} z_{nmq} (t_q - v_q) & \forall n \in N, m \in M \label{eq:APP-1} \\
    s_{n, m+1} & \geq s_{nm} + \sum_{q \in Q} z_{nmq} v_q & \forall n \in N, m \in \{1, ..., \overline{M}-1 \} \label{eq:APP-2} \\
    \sum_{n \in N} \sum_{m \in M} z_{nmq} & = 1 & \forall q \in Q \label{eq:APP-3} \\
    \sum_{q \in Q} z_{nmq} &\leq 1 & \forall n \in N, m \in M \label{eq:APP-4}\\
    \sum_{m \in M} \sum_{q \in Q} (y_{n} - z_{nmq}) &\geq 0 & \forall n \in N \label{eq:APP-5} \\
    z_{nmq} \in \{0, 1\}, & \; y_{nm} \in \{0,1\}, \; s_{nm} \in \mathbb{R} \label{eq:APP-binary}
\end{align}

The objective~(\ref{eq:APP-obj}) minimizes the number of teams necessary. Constraints~(\ref{eq:APP-1}) guarantee that the start moment is after the release time of the corresponding job and before the latest start moment for the corresponding job (i.e. the deadline minus the duration). Constraints~(\ref{eq:APP-2}) enforce that the start moments for one team are sufficiently far apart so that maintenance activities do not overlap. Constraints~(\ref{eq:APP-3}) ensure that every job is assigned to exactly one moment. Constraints~(\ref{eq:APP-4}) make sure that each moment is used for at most one job. Constraints~(\ref{eq:APP-5}) establish that a team can only be used if it is 'active'.  Constraints~(\ref{eq:APP-binary}) ensure that the integer decision variables are also binary.



%% file: 05_MSLCP.tex

\section{Maintenance Scheduling and Location Choice Problem} \label{sec:MSLCP}
The goal of the \textbf{MSLCP} is to find a solution to the \textbf{MLCP} that satisfies predetermined constraints regarding the available number of maintenance teams. To this end, the \textbf{MSLCP} integrates the \textbf{MLCP} and \textbf{APP} in one framework using an approach called Logic-Based Benders' Decomposition (LBBD), which is a generalization of the recognized method called Benders' Decomposition~\citep{hooker2011logic}.. 

\textit{Benders' decomposition} (BD) is a method proposed by \citet{Benders1962} and aims to efficiently solve large-scale linear optimization problems by decomposing the complete problem into a master problem and a sub problem. First, the master problem is solved. Based on the solution of the master problem, a sub problem is identified and solved. Based on the solution of the sub problem, constraints (also called \textit{cuts}) are added to the master problem, which is then solved again. This process continues in an iterative manner. Optimality is reached when the objective value of the master problem is equal to the objective value of the sub problem and the algorithm terminates. In classical BD, cuts are generated via a standard procedure using duality theory. However, in order to do so it requires a specific form for the sub problem (i.e. the complete problem needs to be formulated as one mixed integer program) and it requires that the sub problem be linear and continuous. LBBD does not require that the sub problem take a specific form, at the cost of the fact that it does not have a standard procedure to generate cuts.

Section~\ref{sec:cmslcp_algexplanation} explains the use of LBBD to solve the \textbf{MSLCP}. Section~\ref{sec:capmodel_maths} gives a formal mathematical formulation of the \textbf{MSLCP}. Section~\ref{sec:capmodel_cuts} elaborates on the most important sub procedure of the \textbf{MSLCP}, the cut generation procedure, and presents four different variants for it.

\subsection{\textbf{MSLCP} solution approach}
\label{sec:cmslcp_algexplanation}

The current section proposes an algorithm for the \textbf{MSLCP} using LBBD, in which the \textbf{MLCP} is defined as the master problem and the \textbf{APP} as the sub problem. 

\begin{figure}[H]
    \centering
    \includegraphics{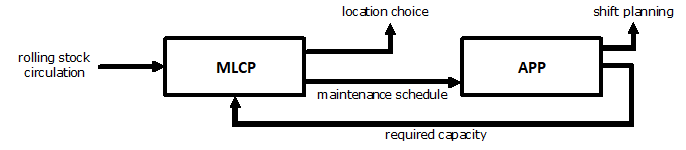}
    \caption{Graphical representation of the \textbf{MSLCP}, demonstrating how it integrates the \textbf{MLCP} and the \textbf{APP}.}
    \label{fig:CMSLCP_flow}
\end{figure}

Figure~\ref{fig:CMSLCP_flow} visualises the cooperation between the \textbf{MLCP} and the \textbf{APP} to include capacity constraints. The maintenance schedule of the \textbf{MLCP} is used to determine the required capacity in the \textbf{APP} model. If the required capacity exceeds the available capacity, the information from the \textbf{APP} is used to add constraints to the \textbf{MLCP} and the \textbf{MLCP} is run again. 

The \textbf{MSLCP} model repeatedly executes the following steps. First,  an empty set of \textit{cuts} is initialized. Second, the \textbf{MLCP} subject to the current set of all generated cuts is solved. Third, by solving APP, a candidate solution for a maintenance schedule, i.e. an assignment of maintenance activities to MOs, is generated.
The algorithm terminates when the \textbf{APP} results in a feasible solution for all time shifts. In that case all constraints in the \textbf{MLCP} and all additional constraints handled by the \textbf{APP} are satisfied and an optimal solution has been determined. Also, it terminates when the current running time exceeds the predetermined maximum running time. Otherwise, it returns to solving \textbf{MLCP} again with newly generated cuts.


\subsection{\textbf{MSLCP} algorithm} \label{sec:capmodel_maths}

\paragraph{Maintenance shifts and maintenance tams} Let $S$ be the set of unique maintenance shifts. Let $\overline{N}$ be the number of teams available at any location. 

\paragraph{Master and sub problem solutions} The \textbf{MSLCP} iteratively solves the master and sub problem. Let $\rho^\kappa$ be the solution of the \textbf{MLCP} after the $\kappa$th iteration of the \textbf{MSLCP} (i.e. this corresponds to a maintenance schedule, which is an assignment of maintenance activities to MOs). Let $Q_{\rho^\kappa}(s)$ be the set of jobs for shift $s \in S$, given the solution of the \textbf{MLCP} $\rho^\kappa$. Let $\mathbf{APP}(Q, \overline{N})$ be the objective value obtained after running the \textbf{APP} for set of jobs $Q$. Use the notation $\mathbf{APP}(Q) = \infty$ if the \textbf{APP} for the set of jobs $Q$ results in an infeasible solution, meaning the required capacity exceeds $\overline{N}$ maintenance teams. To describe the capacity required for a shift $s \in S$, given a solution $\rho^\kappa$ of the master problem, the notation $\mathbf{APP}\left(Q_{\rho^\kappa}(s)\right)$ is used.

\paragraph{Cuts} If $\mathbf{APP}(Q) = \infty$ for a given set $Q$, it can be concluded that the combination of jobs in the set $Q$ results in a violation of the maintenance location capacity. In this case, based on the set $Q$, cuts can be generated according to one of the procedures that are described in Section~\ref{sec:capmodel_cuts}.  A cut indicates a combination of jobs that results in an infeasible solution of the \textbf{APP}. Let $C(Q)$ be the set of cuts based on set $Q$. For any cut $A \in C(Q)$ it holds that $A \subseteq Q$ and $\mathbf{APP}(A) = \infty$. 

Each cut can be translated into a constraint of the \textbf{MLCP} in the following way. Consider a cut $A$. Since $A \subseteq Q$, every element in $A$ signifies a maintenance job which is notated as a tuple $(i, j, K)$ where $i$ is the rolling stock unit, $j$ is the corresponding MO and $K$ is the set of assigned maintenance activities. To include a cut $A$ in the \textbf{MLCP}, the constraint in Equation~(\ref{eq:CMSLCP_cuteq}) needs to be added to prevent the combination of jobs in the cut to show up in a next iteration of the \textbf{MSLCP}. 
\begin{align}
    \sum_{(i, j, K) \in A} \sum_{k \in K} (1-x_{i j k}) \geq 1 \label{eq:CMSLCP_cuteq}
\end{align}

Multiple cuts, for example the set of cuts $C(Q)$, can be added by adding the constraint from Equation~(\ref{eq:CMSLCP_cuteq}) to the \textbf{MLCP} for every cut $A \in C(Q)$.

\paragraph{Iterative procedure} Let $\kappa$ be an index that tracks the current iteration. Let $C^*_{\kappa}$ be the set of cuts generated up to and including the $\kappa$th iteration. Let $C^*_0 = \varnothing$. Let $\ell_0$ be the start time of the algorithm. Let $\ell$ be a parameter restricting the total computation time until the process terminates (if no optimal solution is found earlier).

Pseudo-code for the iterative procedure of the \textbf{CMSCLP} is given in Algorithm~\ref{alg:CMSLCP}.

\begin{algorithm}[H]
\small
	\caption{\textbf{MSLCP} iterative approach} \label{alg:CMSLCP}
	\begin{algorithmic}[1]
		\Function{MSLCP}{$\ell$}
		\State $C^*_0 \leftarrow \varnothing$
		\State $\ell_0\leftarrow$ current time
		\State $\kappa \leftarrow 1$
		\While {$\text{current time } - \ell_0 < \ell$}
		\State compute \textbf{MLCP} solution $\rho^\kappa$, subject to cuts in $C^*_{\kappa-1}$
		\State $C_\kappa^* \leftarrow C_{\kappa-1}^*$
		\For{$s \in S$}
		    \If{$\mathbf{APP}(Q_{\rho^\kappa}(s)) = \infty$}
		    \State $C^*_\kappa \leftarrow C^*_{\kappa} \cup C(J_{\rho^\kappa}(s))$
		    \EndIf
		\EndFor
		\If{$|C^*_{\kappa-1}| = |C^*_{\kappa}|$}
		    \State \textbf{return} $\rho^\kappa$ as the optimal \textbf{MLCP} solution
		\EndIf
		\State $\kappa \leftarrow \kappa+1$
		\EndWhile
		\State \textbf{return} $\rho^\kappa$ as the best found sub-optimal \textbf{MLCP} solution
		\EndFunction

	\end{algorithmic}
\end{algorithm}

The algorithm starts by initializing $C_0^*$, $\ell_0$ and $\kappa$, after which the iterative loop starts. This loop first computes a solution to the \textbf{MLCP} subject to all cuts generated so far. Then, for each shift $s \in S$ in which the required capacity exceeds the available capacity, cuts are generated. The process terminates if either an optimal \textbf{MLCP} solution is found, satisfying all constraints, or if the user-defined maximum running time is exceeded.


\subsection{Cut generation} \label{sec:capmodel_cuts}
If a solution to the \textbf{MLCP} is found that violates the maintenance location capacity constraints, cuts are added to the \textbf{MSCLP} in order to constrain the solution space and prevent such a solution from showing up again. A cut is a set of jobs that cannot occur together since it would result in a violation of available capacity. Cuts result in a restriction of the solution space of the master problem. For a quick convergence of the algorithm, it is desirable to add cuts that are as restrictive as possible. In general, cuts with a smaller amounts of jobs are more restrictive than cuts with larger amounts of jobs. As an example, suppose that the set of maintenance jobs $\{A, B, C \}$ results in an infeasible solution but that the set of maintenance jobs $\{A, B\}$ results in an infeasible solution as well. Both sets of jobs would constitute a valid cut, but the latter set of jobs is smaller, hence more restrictive and as a result more efficient to add.

The remainder of this section proposes four different cut generation procedures: a naive one, then two heuristic ones (a basic and a binary search one), and lastly, a min-cut, which is a more complex one that uses the structure of the problem.

\subsubsection{Naive cut generation} \label{CMSLCP_naivecuts}
Let $Q$ be a set of jobs for that results in a capacity violation, i.e. $\mathbf{APP}(Q) = \infty$. Let $C(Q)$ be the set of cuts generated for this set of jobs. 

Since $Q$ results in an infeasible solution to the \textbf{APP}, this set itself can be added as a cut. Hence, $C(Q) = \{Q \}$.

\subsubsection{Basic Heuristic cut generation} \label{sec:CMSLCPmodel_basicheuristiccut}
To generate smaller cuts compared to the naive procedure, the \textit{Basic Heuristic} cut generation procedure is proposed. This procedure starts with an empty set $\tilde{Q}$ and then moves random jobs iteratively from $Q$ to $\tilde{Q}$. It checks whether the current set of jobs $\tilde{Q}$ results in a feasible solution of the \textbf{APP}. If it does, the current set $\tilde{Q}$ is not yet an appropriate cut since the combination of jobs currently in $\tilde{Q}$ is not infeasible: hence, another job is added in a new iteration. If, on the other hand, it does not, then the current set of jobs is added as a cut to the \textbf{MLCP}.
Pseudo-code for this procedure is presented in Algorithm~\ref{alg:CMSLCP_heuristic_cuts}.

\begin{algorithm}[h]
\small
	\caption{Basic Heuristic cut generation} \label{alg:CMSLCP_heuristic_cuts}
	\begin{algorithmic}[1]
		\Function{Heuristic Cut Generation}{$Q$}
		    \State $\tilde{Q} \leftarrow \varnothing$
		    \While{$\mathbf{APP}(\tilde{Q}) < \infty$}
		    \State pick random $q \in Q$
		    \State $Q \leftarrow Q \setminus q$
		    \State $\tilde{Q} \leftarrow \tilde{Q} \cup q$
		    \EndWhile
		\EndFunction
	\end{algorithmic}
\end{algorithm}

The proposed procedure is guaranteed to terminate since at some point, all jobs from $Q$ are moved to $\tilde{Q}$, meaning that the contents of $\tilde{Q}$ are equal to the initial contents of $Q$. For this set, it is already known that $\mathbf{APP}(Q) = \infty$ since this was required at the start.

The heuristic cut generation procedure can be run multiple times to generate multiple cuts. In general, these cuts are not identical due to the fact that the choice on which job $q \in Q$ to move from $Q$ to $\tilde{Q}$ is random.

\subsubsection{Binary Search Heuristic cut generation} \label{sec:CMSLCPmode_bsheuristiccut} 
The Binary Search Heuristic cut generation procedure uses the same idea as the Basic Heuristic cut generation procedure, but improves upon the efficiency of the former by applying a procedure that is inspired by the principle of binary search (see for example \citet[p.799]{cormen2009introduction}). 

Let $A$ be an initially empty set such that at any moment in the procedure, the jobs in $A$ result in a feasible solution, i.e. $\textbf{APP}(A) < \infty$. Let $B$ be a set of candidate jobs that, when added to the jobs in $A$, at any moment in the procedure results in an infeasible solution: $\textbf{APP}(A \cup B) = \infty$. The algorithm repeatedly splits $B$ into two halves, a left half $B_L$ and a right half $B_R$, and it computes $\textbf{APP}(A \cup B_L)$. If this results in an infeasible solution, i.e. $\textbf{APP}(A \cup B_L) = \infty$, then the set $B_R$ is discarded. In the subsequent iteration of the algorithm the set B of candidate jobs is reduced to $B_L$. If this results in a feasible solution, i.e. $\textbf{APP}(A \cup B_L) < \infty$, some jobs from $B_R$ still need to be added to achieve a 'just infeasible' solution. In this case, the jobs in $B_L$ are all included in the set $A$, and the remaining candidate jobs $B$ to decide on are the jobs $B_R$. The algorithm terminates when $|B| = 1$.
Pseudo-code for the described procedure is given in Algorithm~\ref{alg:CMSLCP_heuristic_cuts_binsearch}.

\begin{algorithm}[h]
\small
	\caption{Binary Search Heuristic cut generation} \label{alg:CMSLCP_heuristic_cuts_binsearch} 
	\begin{algorithmic}[1]
		\Function{Heuristic Cut Generation}{$Q$}
		    \State $A \leftarrow \varnothing$
		    \State $B \leftarrow Q$
		    \While{$|B| > 1$}
		        \State $B_L \leftarrow \varnothing$
		        \State $h \leftarrow \left \lceil \frac{1}{2} |B| \right \rceil$
		        \For{$i \leftarrow 1$ to $h$}
		        \State pick random $j \in B$
		        \State $B_L \leftarrow B_L \cup \{j\}$
		        \State $B \leftarrow B \setminus \{j\}$
		        \EndFor
		        \State $B_R \leftarrow B$
		        \If{ $\textbf{APP}(A \cup B_L) = \infty$}
		        \State $B \leftarrow B_L$
		        \Else
		        \State $A \leftarrow A \cup B_L$
		        \State $B \leftarrow B_R$
		        \EndIf
		    \EndWhile
		    \State \Return $A \cup B$
		\EndFunction
	\end{algorithmic}
\end{algorithm}

The following loop invariants hold (i.e. those expressions are true at the start and end of each iteration):
\begin{itemize}
    \item $\textbf{APP}(A) < \infty$, meaning that the set of jobs in $A$ is feasible
    \item $\textbf{APP}(A \cup B) = \infty$, meaning that when the set of jobs in $B$ is added to the set of $A$, the resulting set of jobs is infeasible.
\end{itemize}

\subsubsection{Min-cut cut generation} \label{sec:CMSLCPmodel_mincut}
In order to find more efficient cuts, the current section designs a procedure that aims to find cuts with a small amount of jobs, by making use of the specific structure of the problem. To this end, the \textit{Relaxed Activity Planning Problem} (\textbf{RAPP}) is defined, which is a relaxation of the \textbf{APP}. In this research, the \textbf{RAPP} is developed for one maintenance team only, although it is expected that the approach can be generalized
to multiple teams.

The benefit of the definition of the \textbf{RAPP} lies in the fact that any infeasible solution to the \textbf{RAPP} is also an infeasible solution to the \textbf{APP}. Recall that cuts need to be generated if the \textbf{APP} is infeasible (see Algorithm~\ref{alg:CMSLCP}). To generate cuts according to the min-cut cut generation procedure, the \textbf{RAPP} is solved. If the \textbf{RAPP} turns out to be infeasible, the min-cut cut generation procedure described in this section can be used. If the \textbf{RAPP} turns out to be feasible, the min-cut cut generation procedure cannot be used and one needs to resort to other cut generation procedures.

The \textbf{RAPP} is a relaxation of the \textbf{APP} in two ways.  First, the \textbf{RAPP} discretizes the planning horizon to a set of \textit{instants}, which are integer minutes, meaning that jobs can only start and end on integer minutes  and job durations should be specified as integers. In the practical context of the railway industry, this is not expected to be problematic since rolling stock units are usually planned per minute. Second, the \textbf{RAPP} allows for preemption of jobs. This means that, unlike in the \textbf{APP}, the work on a job does not need to be performed uninterruptedly.

\paragraph{\textbf{RAPP} definition} The \textbf{RAPP} attempts to assign jobs to as many distinct instants as its duration. This problem can be viewed as a variant of the bipartite matching problem \citep[p. 732]{cormen2009introduction}, where jobs need to be matched to instants, with this difference that jobs in the current problem usually need to be matched to multiple instants instead of only one. The bipartite matching problem is often modeled as a maximum flow problem \citep{fordfulkerson1956}, for which efficient solution algorithms exist \citep[p. 732-735]{cormen2009introduction}. Following this approach, the current research defines the \textbf{RAPP} as a maximum flow problem.

Let $Q$ be the set of jobs, and let $r_q, t_q, v_q$ be the release time, deadline time and duration for job $q \in Q$, respectively, defined in minutes. It is assumed that the duration $v_q$ is integer. Let $P_q$ be the set of instants at which job $q$ is available. This comprises all minutes between $r_q$ and $t_q$ and can be expressed as follows: $P_q = \left\{ x \in \mathbb{N} : \left\lfloor r_q \right\rfloor \leq x \leq \left\lceil t_q \right\rceil  \right\}$. 
Let $P$ be the set of all time instants at which at least one job is available, $P = \cup_{j \in J} P_j$. 



Observe that the \textbf{RAPP} uses discrete time moments (in full minutes) instead of real-valued time moments. Since, in the railway industry, the release and deadline times are usually given in minutes, this is not restrictive.

\paragraph{Step 1: find the maximum flow} Define a source $s$ and a sink $t$ and let $E_G$ be a set of directed edges with capacity $c_e$ for edge $e \in E_G$. Let $G = (N_G, E_G)$ be a directed flow graph, where its set of nodes $N_G$ is defined by $N_G = \{ s \cup Q \cup P \cup t \}$ and its set of directed edges $E_G$ is constructed as follows:
\begin{itemize}
    \item A directed edge $e \in E_G$ from node $s$ to node $j$ for all $q \in Q$ with capacity $c_e = v_q$
    \item A directed edge $e \in E_G$ from node $q$ to $p$ for all $q \in Q$ and $p \in P_q$,  with unit capacity $c_e = 1$. This implies that, for each job, there is a directed edge to each instant at which it is available.
    \item A directed edge $e \in E_G$ from $p$ to $t$ for all $p \in P$, with unit capacity $c_e = 1$.
\end{itemize}

\begin{figure}[H]
    \centering
        \includegraphics[scale=0.5]{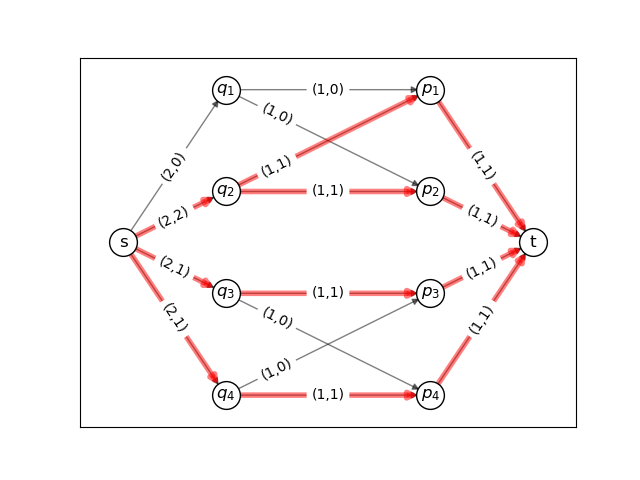}
    \caption{Flow graph $G$ corresponding to the \textbf{RAPP} model. Edges $e$ are annotated $(c_e, f_e)$: the first index represents the edge capacity and the second index represents the assigned edge flow. Red-colored edges (color: see online) represent edges through which a strictly positive flow is assigned.} \label{fig:intelligentcutgen_flowgraphs}
\end{figure}

To illustrate the \textbf{RAPP} cut generation, an example instance is presented where one maintenance team has to perform four jobs: $Q = \{q_1, q_2, q_3, q_4\}$. Jobs $q_1$ and $q_2$ can both be performed at instants $p_1$ and $p_2$ (i.e. $P_1 = P_2 = \{q_1, q_2\}$ and jobs $q_3$ and $q_4$ can be performed at instants $p_3$ and $p_4$ (i.e. $P_3 = P_4 = \{ q_3, q_4\}$). As a result, the set of all instants $P = \{p_1, p_2, p_3, p_4\}$. All jobs have a duration of 2 instants. Figure~\ref{fig:intelligentcutgen_flowgraphs} pictures the associated flow graph.

Once the flow graph has been determined, determine the maximum flow through the flow graph $G$ from the source $s$ to the sink $t$ and denote the resulting flow through each edge $e \in E_G$ by $f_e$. The \textbf{RAPP} is considered to be feasible if and only if the value of the maximum flow equals the sum of all durations, or, equivalently, equals the sum of all capacities on edges departing from $s$, i.e. if and only if 
\begin{align}
    \sum_{e \in E_G} f_e = \sum_{q \in Q} v_q = \sum_{e \in \{(s,v) \in E_G : v \in Q \}} c_e. \label{eq:cmslcp_intelligentcutgen_feasibilitycondition}
\end{align}
The satisfaction of the aforementioned condition(s) represents the fact that all jobs have been completely scheduled. 

Figure~\ref{fig:intelligentcutgen_flowgraphs} pictures the assigned flow on each of the edges in $G$. The maximum flow is 4, whereas the sum of all job durations is 8, meaning that by Equation~(\ref{eq:cmslcp_intelligentcutgen_feasibilitycondition}) the \textbf{RAPP} is not feasible. The remainder of the current section discusses how this infeasible solution can be used to generate cuts.

\paragraph{Step 2: determine the residual graph}
To  find jobs that cannot occur together, the concept of \textit{minimum cuts} from graph theory is used. The capacity of the minimum cut is equal to the value of the maximum flow, and the cut itself provides information about the edges that form a bottleneck in the current graph \citep[p. 269]{taha2011operations}. 

To determine the minimum cut, the concept of \textit{residual graph} is used \citep[p. 716]{cormen2009introduction}. It offers information on how the flow between edges can be changed and represents the amount of possible additional flow through each edge. It may also contain so-called \textit{reverse edges}, that represent the possibility of canceling already assigned flow.

To formally define the concept of the residual graph, let $R$ be a directed graph with the same nodes as $G$ and let its set of edges be denoted by $E_R$, that is, $R = (N_G, E_R)$. Then, the set of edges $E_R$ is constructed as follows. For every edge $e \equiv (u,v) \in E_G$:
\begin{itemize}
    \item there is an edge $e' \equiv (u,v) \in E_R$ with capacity $c_{e'} = c_e - f_e$ if and only if $c_e - f_e > 0$; and
    \item there is an edge $e'' \equiv (v, u) \in E_R$ with capacity $c_{e''} = f_e$ if and only if $f_e > 0$.
\end{itemize}

The nodes that are reachable from $s$ comprise the minimum cut, and the edges connecting one of these nodes to one of the unreachable ones constitute together the bottleneck.

\begin{figure}[H]
    \centering
        \includegraphics[scale=0.5]{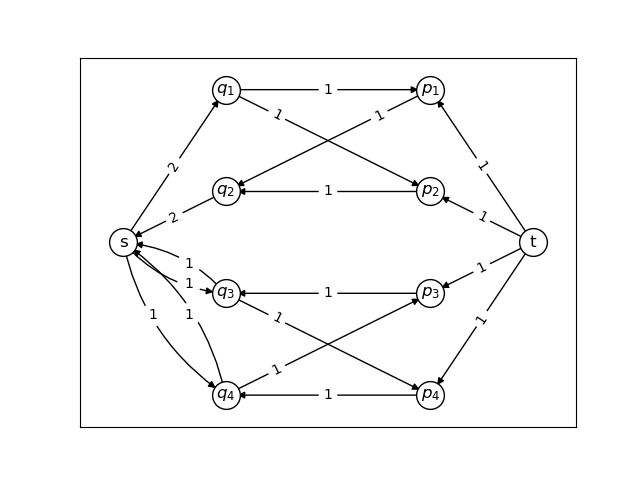} 
    \caption{Residual graph $R$ corresponding to the infeasible solution from Figure~\ref{fig:intelligentcutgen_flowgraphs}. Each directed edge represents the residual capacity between two nodes, if positive.} \label{fig:intelligentcutgen_residualgraph}
\end{figure}

Figure~\ref{fig:intelligentcutgen_residualgraph} displays the residual graph $R$ corresponding to the earlier example in Figure~\ref{fig:intelligentcutgen_flowgraphs}. Take, for instance the positive residual capacity of 2 from $s$ to $q_1$: this signifies that an additional flow can be assigned from $s$ to $q_1$ (corresponding to the situation in which $q_1$ is scheduled). However, in this case, the flow must continue to $p_1$ and $p_2$ (meaning that $q_1$ is scheduled during $p_1$ and $p_2$). This can only be achieved if already assigned flow to $p_1$ and $p_2$ flows back to $q_2$ (signifying that $q_2$, which was formerly scheduled at $p_1$ and $p_2$, is not scheduled anymore) and from there flow further back to the source $s$. The fact that there apparently exists a path from $s$ via $q_1$, $p_1$ and $q_2$ back to $s$ is an important observation: it signifies that $q_1$ and $q_2$ are conflicting. This, in turn, means that $q_1$ and $q_2$ cannot be scheduled together and can be added as a cut. In fact, all jobs on every path starting from $s$ and returning to $s$ constitute an infeasible combination of jobs.

\paragraph{Step 3: define the Reachable Components graph} To formalize the idea of conflicting jobs, the \textit{Reachable Components graph} $H$ is introduced. Its aim is to separate components that define different combinations of jobs, each of which cannot occur together (i.e. result in an infeasible solution of the \textbf{RAPP}). Let $H$ be a directed graph and let it have the same nodes as $G$ and with the set of edges $E_H$, i.e. $H = (N_G, E_H)$. Let $E_H$ contain all edges in $R$ that are not connected to the source $s$ or sink $t$, that is, $E_H = \{ (u,v) \in R: u \notin \{s, t\}, v \notin \{s, t \} \}$. Let $D(F, n)$ be the set of all nodes reachable in some graph $F$ starting from some node $n$ (also called the \textit{descendants} of $n$ in $F$). This set of reachable nodes can be obtained efficiently by the application of a depth-first search \citep[p. 603-606]{cormen2009introduction}.

From this, a set of cuts can be determined. Note that all separate sets of reachable nodes can be obtained by starting at some job $j \in J$ that is reachable from $s$ in $R$ and obtaining all jobs among its descendants. In other words, for all $q \in Q : (s, q) \in R$ the set $C_q = \{ q \cup (D(H, q) \cap Q) \}$ comprises a set of jobs that cannot occur together. These jobs result in an infeasible \textbf{RAPP} solution and, as a consequence, in an infeasible \textbf{APP} solution; hence, they can be added as a cut for the \textbf{MLCP}.

\begin{figure}
\centering
        \includegraphics[scale=0.5]{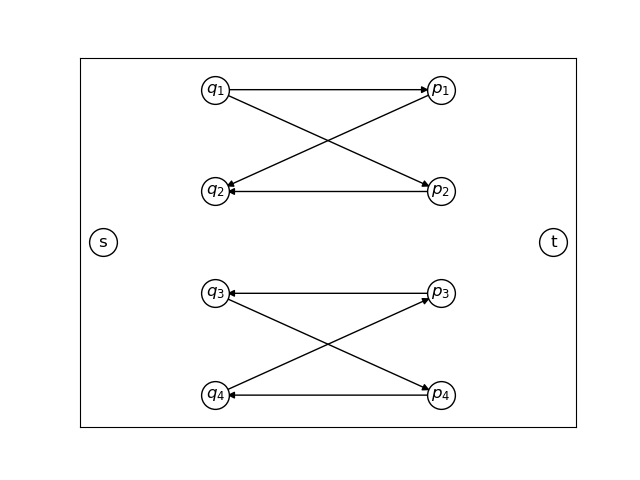}
        \caption{Reachable Components graph $H$, separating the various reachable components that are reachable from $s$.}
        \label{fig:intelligentcutgen_reachablecomponents}
\end{figure}

To demonstrate the process of the determination of these cuts, return once again to the previous example. Figure~\ref{fig:intelligentcutgen_reachablecomponents} presents the graph $H$ with two different components. In $R$, the nodes $q_1, q_3$ and $q_4$ are reachable from $s$. Hence, the cuts generated in this way are $\{q_1, q_2 \}, \{ q_3, q_4\}$ and $\{ q_4, q_3\}$. This shows that $q_1$ and $q_2$ cannot occur together, and similarly that $q_3$ and $q_4$ cannot occur together.

\paragraph{Step 4: cut set post-processing} All cuts according to the above described procedure can be added to the \textbf{MLCP}, but some of these may be superfluous. First, the same cuts may be generated more than once. Second, some cuts may be generated while a more specific cut is also generated: for example, consider the generation of two cuts, the first with jobs $X, Y$ and $Z$ and the second with jobs $X$ and $Y$. The latter makes the former redundant. 

To remove redundant cuts, a straightforward procedure is applied that iteratively adds cuts only if it is not a superset of a more efficient cut that was already added. To this end, let $C$ be the set of all cuts generated by the \textbf{RAPP} and let $\Tilde(C)$ be the set of cuts with all redundant cuts from $C$ removed. Algorithm~\ref{alg:intcuts_removeredundant} gives pseudo-code for this procedure.

\begin{algorithm}[H]
\small
	\caption{Remove redundant cuts after min-cut cut genreation} \label{alg:intcuts_removeredundant}
	\begin{algorithmic}[1]
		\Function{Remove redundant cuts}{$C$}
		\State sort $C$ by the cardinality of all its elements $c \in C$
		\State $\Tilde{C} \leftarrow \varnothing$
		\For{$c \in C$}
		    \State $\text{add} \leftarrow \textbf{true}$
		    \For{$\Tilde{c} \in \Tilde{C}$}
		        \If{$c \supseteq \Tilde{c}$}
		        \State $\text{add} \leftarrow \textbf{false}$
		        \EndIf
		    \EndFor
		    \If{$\text{add} = \textbf{true}$}
		        \State $\Tilde{C} \leftarrow \Tilde{C} \cup \{ c \}$
		    \EndIf
		\EndFor
		\State \textbf{return} $\Tilde{C}$
		\EndFunction
	\end{algorithmic}
\end{algorithm}



%% file: 06_Results.tex

\section{Experimental results} \label{sec:results}
The current section investigates the performance of the \textbf{MSLCP} model. It considers a smaller-scale instance to investigate how well the model is able to find an optimal solution, and it considers a larger-scale (and hence more realistic) instance to investigate how quickly the model is able to converge to a solution that, although it may be sub-optimal, is useful in practice.

\subsection{Scenario set-up} \label{sec:resultsCMLCP_expdesign}

The \textbf{MSLCP} framework is demonstrated on realistic scenarios from the Dutch railways. The problem instance considered in the current section uses a rolling stock circulation originating from the main Dutch railway operator, Netherlands Railways (NS). Specifically, it uses so-called BasisDag update (BDu) data for the period between 10-4-2018 until 16-4-2018. 
In particluar,
4 rolling stock types are considered: ICM4, DDZ4, DDZ6 and DD-AR3. These rolling stock types are chosen in such a way that they result in some maintenance location capacity issues, especially at maintenance location Zwolle (Zl). This comprises a total of 137 rolling stock units. 

The planning horizon is set to 7 days, equal to the total number of days in the input data. The set of nighttime maintenance locations $L^N$ and the set of potential daytime maintenance locations $L^D$ are assumed to be equal to the set of all locations in the BDu. It is assumed that 5 locations can be opened for daytime maintenance at maximum, i.e. $L^D_{max} = 5$. Two maintenance types are included, maintenance type A having a duration of 30 minutes and an interval of 24 hours, maintenance type B having a duration of 60 minutes and an interval of 48 hours. Rolling stock units are assumed to be as-good-as-new at the start of the planning horizon, i.e. $b_{ik} = 0 \text{ for all }  i \in I, k \in K$. The technical parameter $\varepsilon$ has a value of $\varepsilon = 0.001$. It is assumed that at each maintenance shift (at each location, on each day), one maintenance team is available, i.e. $\overline{N}=1$. Note that this assumption is actually necessary for the min-cut cut generation procedure, which is only defined for one maintenance team. Unless stated otherwise, the running time is restricted to $\ell = 2$ hours.

Two scenarios are constructed.
In both scenarios, the capacity of daytime maintenance shifts is considered only, while the capacity of nighttime maintenance shifts is ignored. This choice is reasonable in the light of the gradual introduction of a policy of daytime maintenance in practice, where capacity for daytime maintenance at first is limited. The set of shifts $S$ is dependent on the scenario used and is discussed below.

First, the \textit{single-shift scenario} aims to investigate how quickly the \textbf{MSLCP} converges to optimality, that is a solution without capacity constraint violations, for the proposed four types of cut generation procedures. Reaching the optimal solution may take long, and since the time to find such an optimal solution relates to the number of maintenance shifts for which capacity constraints are imposed, this scenario focuses at one maintenance shift:  the daytime maintenance shift at maintenance location Zl on 11-4-2018. The set of shifts $S$ contains only this maintenance shift. This particular shift was selected after exploratory experiments showed that determination of the required capacity using the \textbf{APP} took most time for this shift, and that it contains relatively many maintenance activities that are also overlapping. As a result, this maintenance shift shows to be 'hard' to solve. Using a hard-to-solve maintenance shift in the single-shift scenario enables to investigate the performance of various cut generation procedures in solving a capacity violation of a specific shift, as opposed to using a maintenance shift that would be more easy to solve which would make it harder to demonstrate differences between cut generation procedures. In the single-shift scenario, 10 different \textit{cut generation variants} are investigated: cut generation by the naive cut generation procedure (one variant), by the Basic Heuristic cut generation procedure, for 1, 2, 5 and 15 cuts (four variants), by the Binary Search Heuristic cut generation procedure, for 1, 2, 5 and 15 cuts (four variants) and by the min-cut cut generation procedure (one variant). The most important performance indicator is the convergence of the current objective value of the \textbf{MLCP} as a function of either computation time or the number of iterations.

Second, the \textit{all-shifts scenario} accounts for the fact that, in practice, reaching an optimal solution may take too much computation time. 
In particular, a sub-optimal solution with limited capacity violations obtained quickly may, in practice, be preferred over an optimal solution without any violations taking excessive computation time. Therefore, it is worthwhile to know how quickly the number of capacity violations can be reduced.
To investigate the performance of the proposed algorithms in this practical case, the second scenario presents a more realistic setup and includes all daytime maintenance shifts.
The set of shifts $S$ contains maintenance shifts for all possible combinations of maintenance location and date in the planning horizon.  In the all-shifts scenario, three cut generation variants are considered: the naive cut generation procedure, the Binary Search Heuristic cut generation procedure for 15 cuts, and the min-cut cut generation procedure. Compared to the one-shift scenario, the seven other heuristic cut generation methods, i.e. Basic Search Heurstic procedure with 1, 2, 5 and 15 cuts and the Binary Search Heuristic procedure with 1, 2 and 5 cuts, are left out. The Binary Search Heuristic cut generation procedure with 15 cuts outperformed the other seven procedures, therefore only this one is considered.
%
The most important performance indicator is the number of shifts for which the required capacity exceeds the available capacity, i.e. the number of capacity violations, as a function of either time or the number of iterations. It has been verified that a solution to the \textbf{MLCP} can be obtained in which all maintenance is performed during nighttime is feasible. Hence, if one would allow for enough computation time, the number of capacity violations would converge to zero with certainty.

The \textbf{MSLCP}, \textbf{MLCP}, \textbf{APP} and \textbf{RAPP} are implemented implemented using Python and solved using Gurobi. For the implementation of the \textbf{RAPP}, the package NetworkX \citep{hagberg2008exploring} is used. The corresponding maximum flow problem is solved using the preflow-push algorithm (see e.g. \citet[p. 765]{cormen2009introduction}), that is included in the implementation of NetworkX. An implementation code of the \textbf{MSLCP} model is provided by \citet{GithubCMLCP2020}. For reasons of confidentiality, the actual data could not be provided, but synthetic data is made available instead.

\subsection{Results} \label{sec:resultsCMLCP_results}
This section provides the results generated for the \textbf{MSLCP}. Section~\ref{sec:results_workings} provides an illustrative example for a particular maintenance shift, demonstrating how the \textbf{MSLCP} (described in Algorithm~\ref{alg:CMSLCP}) is able to find a schedule that satisfies capacity constraints. 
Then, Section~\ref{sec:CMLCPresults_results_singleshift} and 
Section~\ref{sec:CMLCPresults_results_allshifts} give results for the two scenarios proposed in 
Secton~\ref{sec:resultsCMLCP_expdesign}.

\subsubsection{Illustrative schedule} \label{sec:results_workings}
To demonstrate the workings of the \textbf{MSLCP}, the schedule for the daytime maintenance shift on 13-4-2018 at maintenance location Zl is examined. This shift is useful for demonstration purposes since it consists of an interesting variety of maintenance jobs. 

\begin{figure}
    \centering
    \begin{subfigure}{.48\textwidth}
        \includegraphics[scale=0.45]{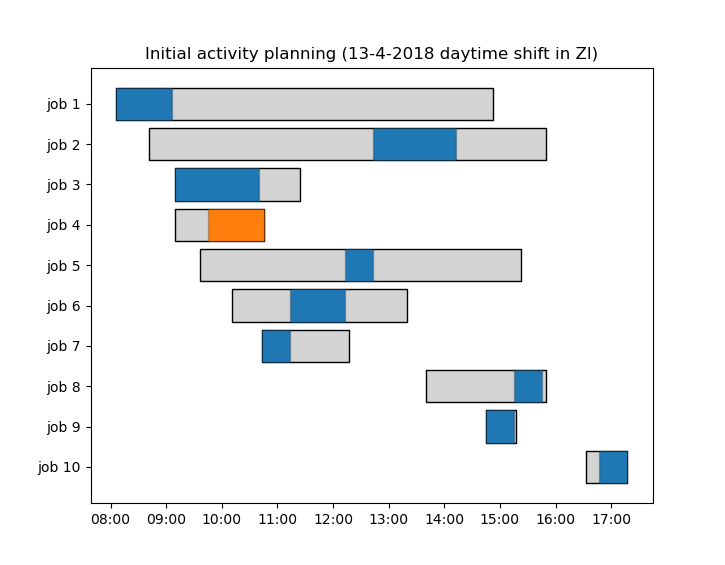} 
    \end{subfigure}
    \begin{subfigure}{.48\textwidth}
        \includegraphics[scale=0.45]{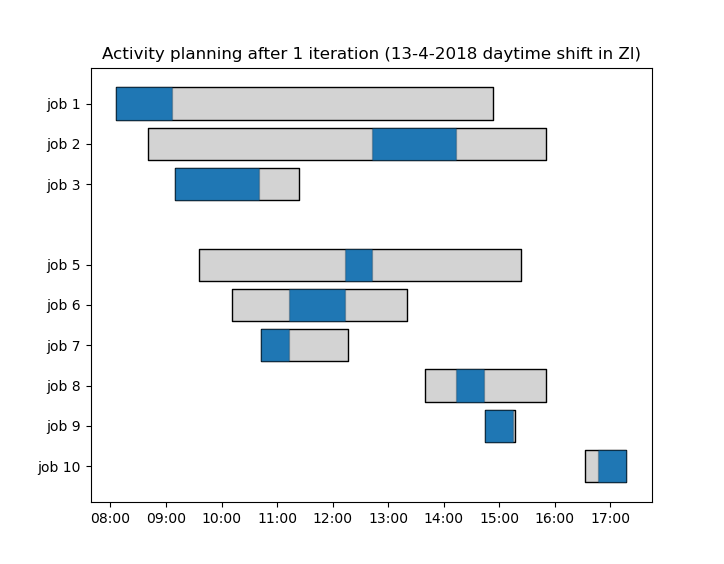}
    \end{subfigure}
    \caption{The maintenance schedule for the daytime maintenace shift in Zl on 13-4-2018 in the initial situation (left) and after one iteration of the \textbf{MSLCP} (right). The \textbf{MSLCP} finds a solution that requires only one maintenance team after one iteration. The x-axis shows the time period 8:00-17:00, and y-axis the jobs for this period. Grey boxes represent the intervals of the maintenance opportunity during which a maintenance job should take place. Blue represents the actual scheduled time if the job is performed by team 1 and orange if it is performed by team 2.} \label{fig:resultsMSLCP_illustrativeschedule}
\end{figure}

Figure~\ref{fig:resultsMSLCP_illustrativeschedule} demonstrates schedules for this particular maintenance shift, computed by the \textbf{APP} model. The maintenance jobs assigned to this shift are determined by the initial \textbf{MSLCP} solution (Figure~\ref{fig:resultsMSLCP_illustrativeschedule}, left) and after one iteration (Figure~\ref{fig:resultsMSLCP_illustrativeschedule}, left). 

The length of a maintenance job varies: it may represent a type A maintenance activity (30 minutes), a type B maintenance activity (60 minutes) or a combination of both (90 minutes). The interval of the maintenance opportunity (MO) during which a maintenance job should take place is depicted in grey, and the actual scheduled time is shown in blue if it is performed by team 1 and in orange if it is performed by team 2. 

In the initial schedule of the maintenance shift (Figure~\ref{fig:resultsMSLCP_illustrativeschedule}, left),  
observe that the combination of maintenance jobs that need to be performed cannot be fulfilled by only one team and requires two teams instead (see the orange job between 9:45 and 10:45). Since the required capacity exceeds the available capacity of maintenance teams, the \textbf{MSLCP} procedure generates cuts and finds a new \textbf{MLCP} solution. 
In the new \textbf{MLCP} solution (Figure~\ref{fig:resultsMSLCP_illustrativeschedule}, right), the assignment of maintenance activities to MOs has changed in such a way that maintenance job 4 is not part of this maintenance shift anymore. As a result, a schedule can be created that requires only one maintenance team and the \textbf{MSLCP} procedure has solved a capacity violation.

\subsubsection{Single-shift scenario} \label{sec:CMLCPresults_results_singleshift}
\begin{figure}
    \centering
    \begin{subfigure}{.48\textwidth}
        \includegraphics[scale=0.45]{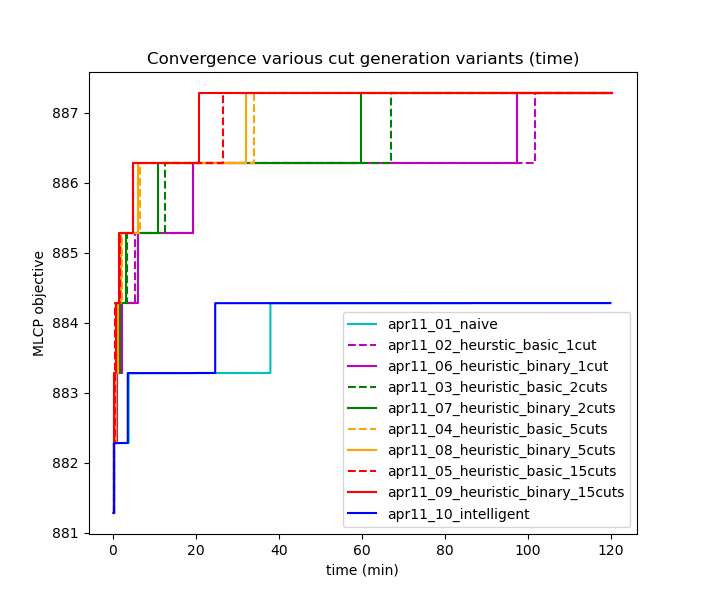} 
    \end{subfigure}
    \begin{subfigure}{.48\textwidth}
        \includegraphics[scale=0.45]{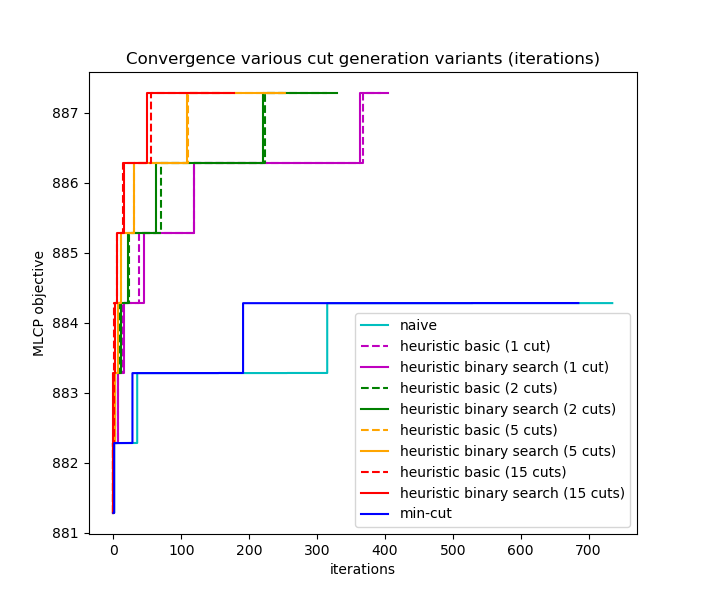}
    \end{subfigure}
    \caption{Convergence of the \textbf{MSLCP} in the single-shift scenario. For each cut generation variant, the course of the value of the \textbf{MLCP} is displayed as a function of elapsed time (left) and as a function of the current iteration (right).} \label{fig:resultsCMLCP_singleshift_conv}
\end{figure}

Figure~\ref{fig:resultsCMLCP_singleshift_conv} provides a graphical representation of the development of
the LBBD framework for \textbf{MSLCP}, i.e.
the objective function of 
the latest \textbf{MLCP} solution over time and over multiple iterations. It shows that all heuristic cut generation variants achieved an objective of approximately 887. To be precise, the heuristic cut generation variants' final objective values are between 887.277 and 887.281\footnote{Recall that the objective value is mainly composed of the total number of daytime activities. The reason that the value is nonetheless not integer is due to the fact that, besides a unit value for each daytime activity, also a value $\epsilon$ is added for every performed maintenance activity (see Section~\ref{sec:MLCP}).}, thereby coming closest to the (unknown) optimal value and providing a lower bound (887.281) for it. Of these heuristic cut generation variants, the variants with higher number of cuts reach this objective value 
faster (i.e. in less time and in less iterations) than the variants with lower number of cuts. The cut generation variant Binary Search Heuristic with 15 cuts performed the best, i.e. it reached the value of 887 within the least amount of time (22 minutes) and in the least amount of iterations (50 iterations). All  heuristic cut generation procedures reached the objective value of 887 within two hours, unlike the min-cut cut generation procedure (which reached the objective value of 884 within 25 minutes in 192 iterations), and the naive cut generation procedure (which reached the objective value of 884 within 39 minutes in 316 iterations).

When comparing the Binary Search Heuristic with the Basic Heuristic, it is found that their convergence is similar in terms of iterations, but that the convergence of the Binary Search Heuristic is a bit quicker time-wise. This is an indication that the improvement per iteration is comparable for both, but that the time consumed per iteration is less for the Binary Search Heuristic. 

As indicated, the single-shift scenario focuses at one particular maintenance shift and attempts to find a solution in which the capacity constraint for this shift is met. For this goal, the heuristic cut generation variants outperform both the naive and the min-cut cut generation variants, in time as well as in number of iterations. For the latter two, however, much more iterations were performed. This is an indication that, despite the fact that the computation time per iteration is lower, the cuts produced in each iteration by the min-cut and naive cut generation procedures contribute to a lesser extent to the convergence of the \textbf{MSLCP} than in the heuristic cut generation procedures.

In an attempt to 
find an optimal value to benchmark the cut generation variants, the best-performing cut generation variant (heuristic binary search with 15 cuts) was run for 14 hours. Still, no optimal solution to the \textbf{MSLCP} was found, although this run did provide a new lower bound to the optimal objective value of 888.279.

Table~\ref{tab:CMLCPresults_singleshift_comptime} displays the computation time per iteration in various sub processes. 

\begin{table}[H]
\centering
\small
\begin{tabular}{llllll} \FL
& \multicolumn{4}{c}{sub processes} & \\
\cline{2-5}
                                     & \textbf{MLCP} & \textbf{APP} & cut gen. & other & total \ML
naive                     & 9.3      & 0.3    & 0.0          & 0.1   & 9.8   \\
Basic Heuristic (1 cut)     & 16.6     & 0.5    & 1.0          & 0.2   & 18.2  \\
Basic Heuristic (2 cuts)     & 20.2     & 0.5    & 2.0          & 0.2   & 22.9  \\
Basic Heuristic (5 cuts)     & 22.0     & 0.5    & 5.5          & 0.2   & 28.3  \\
Basic Heuristic (15 cuts)    & 26.2     & 0.5    & 17.4         & 0.2   & 44.3  \\
Binary Search Heuristic (1 cut)  & 16.3     & 0.5    & 0.9          & 0.2   & 17.8  \\
Binary Search Heuristic (2 cuts)  & 19.4     & 0.5    & 1.9          & 0.2   & 21.9  \\
Binary Search Heuristic (5 cuts)  & 23.1     & 0.6    & 5.0          & 0.2   & 28.8  \\
Binary Search Heuristic (15 cuts) & 25.0     & 0.5    & 14.9         & 0.1   & 40.5  \\
Binary Search Heuristic (15 cuts) extended* &	87.8 & 	0.7 & 17.1 & 0.2 & 105.7 \\
min-cut            & 9.5      & 0.4    & 0.5          & 0.2   & 10.5 \LL
\end{tabular}
\caption{Computation time per iteration for each cut generation variant, in seconds, decomposed into the main contributing processes to the computation time: the computation of an \textbf{MLCP} solution subject to all cuts generated, the determination of a capacity violation using the \textbf{APP}, and the cut generation process itself, and other processes. The last relates to all remaining computations, such as results storage. 
*The result for the extended computation time of 14 hours.} \label{tab:CMLCPresults_singleshift_comptime}
\end{table}

It can be observed that the naive and min-cut cut generation variants require the least time per iteration. This is in correspondence with the fact that in these variants many iterations could be run within 2 hours (see Figure~\ref{fig:resultsCMLCP_singleshift_conv}). 

Moreover, Table~\ref{tab:CMLCPresults_singleshift_comptime} shows that the generation of cuts in the Basic Heuristic version requires somewhat more time than the Binary Search Heuristic. This concurs with the expectation that can be drawn from the design of both heuristics: the Binary Search Heuristic improves upon the Basic Heuristic in the sense that it requires less iterations to generate a cut.  Also, the iterations of the heuristic cut generation variants take more time for higher numbers of cuts, which is a direct result of the time it takes to generate more cuts.

The average running time of the \textbf{APP}, necessary to determine whether capacity of a maintenance shift is violated, is  well below one second consistently over all cut generation variants. 

\begin{figure}
    \centering
        \includegraphics[scale=0.5]{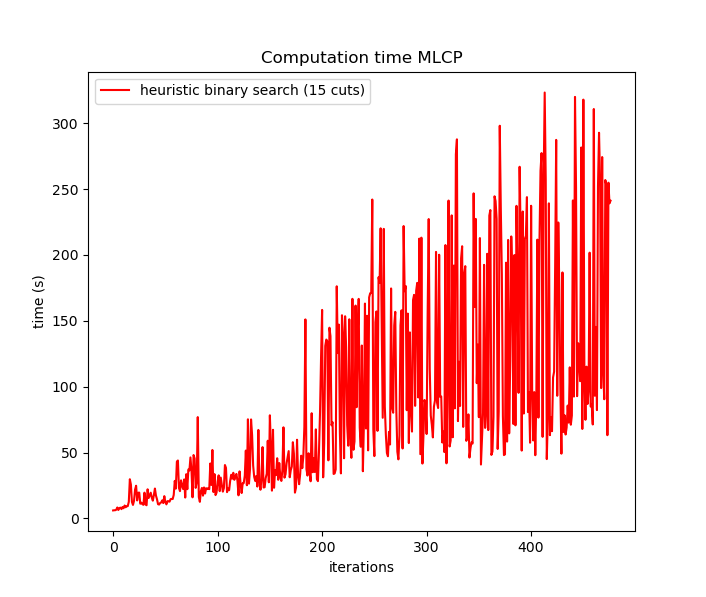}
        \caption{Computation time of the \textbf{MLCP} in seconds, per iteration of the \textbf{MSLCP} for the Binary Search Heuristic cut generation version with 15 cuts per iteration, in an extended run of 14 hours.}         \label{fig:resultsCMLCP_singleshift_MLCPrunningtime}
\end{figure}

The most time is consumed by solving the \textbf{MLCP}. Interestingly, the \textbf{MLCP} takes more time to run in the heuristic cut generation variants than it does in the naive and min-cut cut generation variants. To understand this, it is relevant to look at the computation time of the \textbf{MLCP} for the extended run of 14 hours. Figure~\ref{fig:resultsCMLCP_singleshift_MLCPrunningtime} presents it as a function of the current iteration. It shows that the running time of the \textbf{MLCP} (as well as its variance) increases for later iterations. The  explanation for this is that due to the added cuts, the \textbf{MLCP} becomes increasingly constrained and solving it becomes increasingly difficult. This leads to higher computation times for the \textbf{MLCP}. 

\subsubsection{All-shifts scenario} \label{sec:CMLCPresults_results_allshifts}
Section~\ref{sec:CMLCPresults_results_singleshift} examined the results obtained by using the \textbf{MSLCP} in a context where the number of available maintenance teams of only one single maintenance shift was constrained. The capacity constraint for this single maintenance shift appeared to be highly complicating. This means that many iterations of the \textbf{MSLCP} are necessary to reduce the number of required maintenance teams (so that, in turn, this number meets the number of available maintenance teams). However, in realistic cases, it is not the case that the constraints of each maintenance shift are as complicating. Therefore, the current section considers the all-shifts scenario, attempting to solve the capacity violations for all shifts of the problem instance. 

As a reference, the initial \textbf{MLCP} solution of the scenarios is used, i.e. without initial cuts. To  obtain numbers of the required maintenance teams, the APP was are solved for each maintenance shift.
In this solution, there are 34 daytime maintenance shifts to which at least one maintenance activity is assigned, of which 13 require 1 maintenance team, 19 require 2 maintenance teams and 2 require 3 maintenance teams. Thus, given the capacity of one maintenance team per shift, there are 21 shifts with capacity violation.

\begin{figure}[H]
    \centering
    \begin{subfigure}{.48\textwidth}
        \includegraphics[scale=0.45]{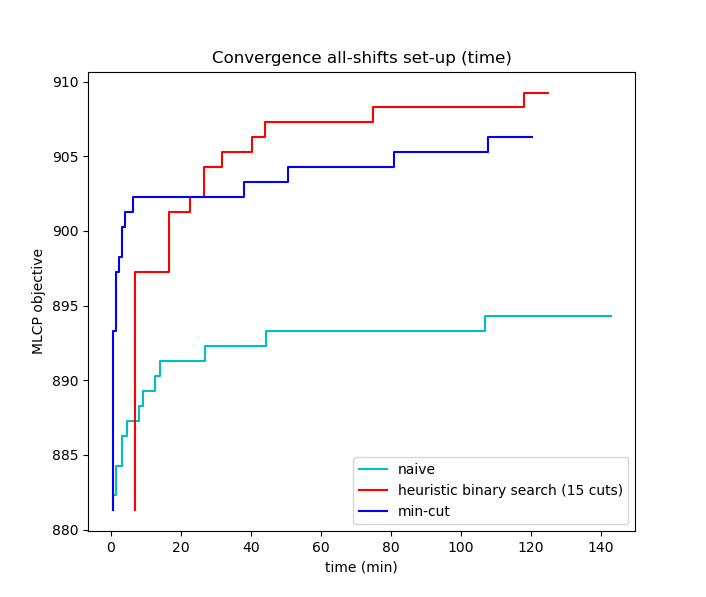}
    \end{subfigure}
    \begin{subfigure}{.48\textwidth}
        \includegraphics[scale=0.45]{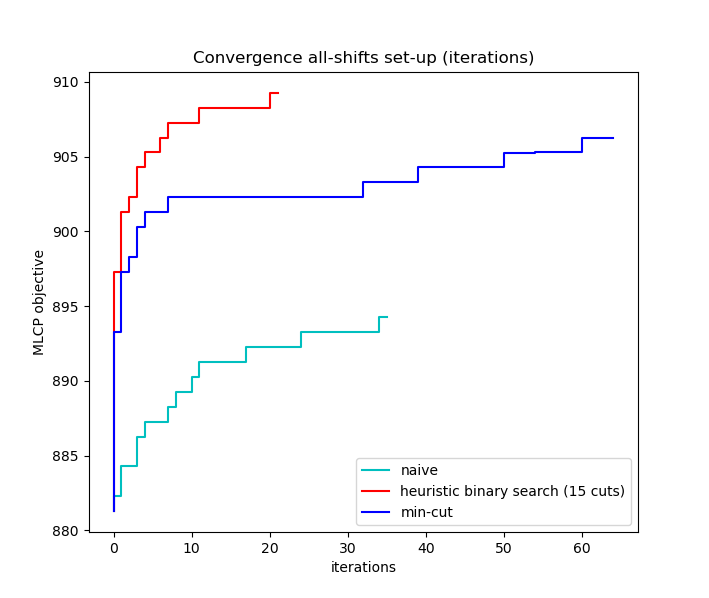}
    \end{subfigure}
        \caption{Convergence of the \textbf{MSLCP} in the all-shifts set-up. For three generation variants, the course of the value of the \textbf{MLCP} is displayed as a function of elapsed time (left) and as a function of the current iteration (right).} \label{fig:resultsCMLCP_allshifts_conv}
\end{figure}

\paragraph{Convergence}
As in the single-shift scenario, no optimal solutions were found within the running time restriction of two hours. Figure~\ref{fig:resultsCMLCP_allshifts_conv} displays the convergence of the \textbf{MLCP} objective value in the all-shifts set-up.

At first, it can be noted that the course of the \textbf{MLCP} objective value for the three investigated cut generation variants is similar as 
for the single-shift scenario. As a result of the added cuts, the value of the objective of the \textbf{MLCP} gradually increases. 
It can be seen that, in terms of iterations, the course of the \textbf{MLCP} objective at the beginning of the run is very similar for the min-cut and Binary Search Heuristic cut generation processes.

In the first couple of iterations (right side of Figure~\ref{fig:resultsCMLCP_allshifts_conv}), both procedures are equally capable of detecting 'simple' infeasible combinations of jobs that, when added as a cut to the \textbf{MLCP}, immediately cause a unit step in the \textbf{MLCP} objective. The min-cut cut generation procedure has an advantage, since its running time per iteration is shorter. This is reflected in the left side of Figure~\ref{fig:resultsCMLCP_allshifts_conv}, where the increase in objective value is quicker in case of the min-cut cut generation procedure. In a later stage, however, the cuts added by the Binary Search Heuristic cut generation procedure yield a better convergence of the \textbf{MLCP} (right side of Figure~\ref{fig:resultsCMLCP_allshifts_conv}). Hence, from a time perspective, in a later stage the Binary Search Heuristic cut generation procedure overtakes the min-cut cut generation procedure (as can be seen in the left of Figure~\ref{fig:resultsCMLCP_allshifts_conv}).

\paragraph{Number of capacity violations}
By comparing the solutions of the \textbf{MSLCP} model (from the last iteration) and the \textbf{MLCP} (i.e. the initial solution), it has been shown that optimality is not reached within 2 hours of computation time, which means that even the last solution obtained after 2 hours of computation time contains maintenance shifts for which the capacity is violated. 

\begin{figure}[H]
    \centering
    \begin{subfigure}{.48\textwidth}
        \includegraphics[scale=0.45]{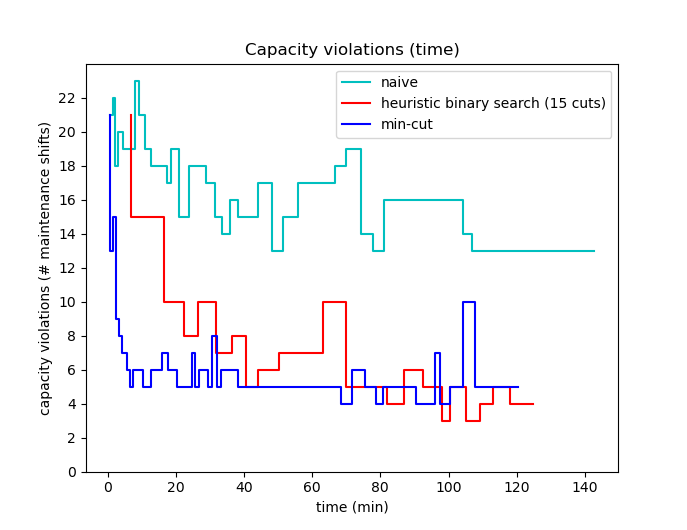}
    \end{subfigure}
    \begin{subfigure}{.48\textwidth}
        \includegraphics[scale=0.45]{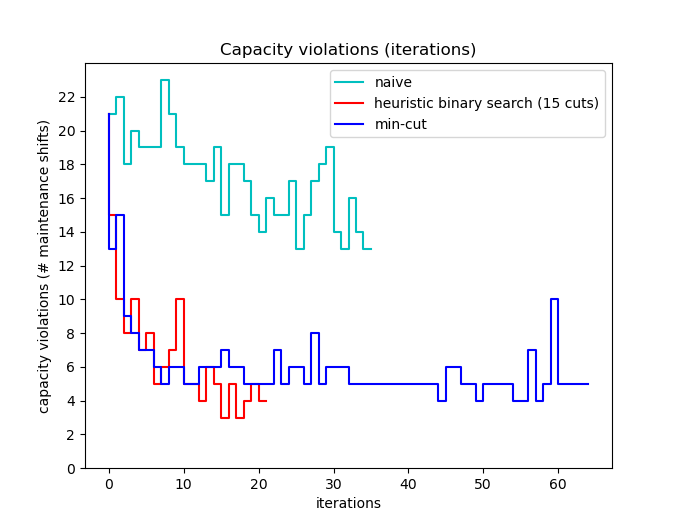}
    \end{subfigure}
    \caption{Number of shifts for which the capacity is violated (i.e the required capacity is more than 1 maintenance team), for three cut generation variants, as a function of elapsed time (left) and as a function of the current iteration (right). The naive cut generation procedure took longer than 2 hours since solving the \textbf{MLCP} in the last iteration (that started before the threshold of 2 hours of running time) took very long; the process terminated as soon as this iteration was finished.} \label{fig:CMLCPresults_allshifts_capviolation}
\end{figure}

A deeper look into the development of the number of capacity violations is taken. Figure~\ref{fig:CMLCPresults_allshifts_capviolation} presents the number of capacity violations as a function of elapsed time and as a function of the current iteration for all cut generation procedures.

First, it becomes clear that, for each cut generation procedure, the number of maintenance shifts for which the capacity is violated starts at 21 (the initial number of capacity violations) and is decreasing. However, the decrease is not strictly monotonic. The added cuts as a result of the violation of capacity in one of the maintenance shifts, may induce a new \textbf{MLCP} solution that assigns maintenance in such a way that the capacity of maintenance shift which was formerly sufficient, now becomes violated.

The naive cut generation variant is clearly the worst performing. After two hours of running time, it contains considerably more maintenance shifts for which capacity is violated than the other two cut generation variants. More strikingly is the development of the number of violations in the min-cut cut generation variant compared to the Binary Search Heuristic cut generation variant.  When looking at the development in terms of the elapsed time, the capacity violations in the min-cut cut generation variant decrease at the beginning much more sharply than in the binary search cut generation variant, after which they in both remain constant for around 5 capacity violations. The min-cut cut generation variant found a solution with 5 maintenance shift violations or less after 7.6 minutes, whereas the Binary Search Heuristic cut generation procedure found such a solution only after approximately 44.2 minutes. The practical implications of this are relevant: when no feasible solution can be obtained in reasonable time, the preferred option is to get a good sub-optimal solution as quick as possible. The min-cut cut generation procedure seems better suited for this goal. 

To gain a little more understanding on this behavior, observe also the capacity violations as a function of the current iteration. At the beginning, the Binary Search Heuristic and min-cut cut generation variants show a similar path. This implies that, in each iteration, the resulting cuts in both variants lead to similar benefits in the reduction of capacity violations. However, the running time of the min-cut cut generation procedure per iteration is considerably lower than in the Binary Search Heuristic cut generation procedure, leading to a better performance in terms of computation time.

In the first couple of iterations, both procedures are equally capable of detecting 'simple' infeasible combinations of jobs that, when added as a cut to the \textbf{MLCP}, immediately cause a unit step in the \textbf{MLCP} objective. The min-cut cut generation procedure has an advantage, since its running time per iteration is shorter. In a later stage, however, the cuts added by the Binary Search Heuristic cut generation procedure yield a better convergence of the \textbf{MSLCP}. Hence, in a later stage the Binary Search Heuristic cut generation procedure overtakes the min-cut cut generation procedure (as can be seen in the left of Figure~\ref{fig:resultsCMLCP_allshifts_conv}). Hence, when an application does not require all capacity violations to be solved, the min-cut cut generation procedure is preferred since it reduces the number of capacity violations most quickly. However, when all capacity violations need to be solved, the binary search heuristic cut generation procedure has a better performance.

%% file: 07_Conclusion.tex

\section{Conclusion} \label{sec:conclusion}

The current work proposes the  Maintenance Scheduling and Location Choice Problem  \textbf{MSLCP}, which provides a rolling stock maintenance schedule and a maintenance location choice, taking into account the available capacity of maintenance locations, measured in the number of available maintenance teams. It combines the \textbf{MLCP} model, introduced in \citet{zomer2020MLCP}, and the \textbf{APP} model (proposed in the the current research) using a framework called Logic-Based Benders' Decomposition. In addition, four different procedures for the generation of \textit{cuts} are proposed.

The \textbf{APP} provides the required number of maintenance teams quickly, i.e. within seconds for realistic problem sizes. This is beneficial since in the \textbf{MSLCP} context it needs to be run for every iteration and therefore contributes to the efficiency of the \textbf{MSLCP} model. In addition, it provides an optimal maintenance shift planning. As such, it is not only a valuable addition to the \textbf{MLCP}, but it can also be useful in operational contexts where a shift planning is required.  

The \textbf{MSLCP} algorithm is designed to find a solution to the \textbf{MLCP} that includes maintenance location capacity. It proposes a framework to incorporate complex constraints in scheduling problems. The current research investigates the performance of the \textbf{MSLCP} in two scenarios.  The first scenario focuses on one particular maintenance shift for which the capacity is violated, demonstrating that the binary search heuristic cut generation procedure with 15 cuts is the most promising procedure to solve the capacity violations of a hard-to-solve instance. The second scenario focuses on solving the capacity issues in multiple maintenance shifts. The number of maintenance shifts for which the required capacity exceeded the available capacity could be reduced from 21 to 5 in less than 8 minutes using the min-cut cut generation procedure. Hence, the min-cut cut generation procedure is able to quickly decrease the number of maintenance shifts necessary to a reasonable amount, but when solving a hard maintenance shift to optimality, the min-cut cut generation procedure is outperformed by the binary search heuristic cut generation procedure with 15 cuts. 


The current research has some limitations. First, it is assumed that maintenance jobs need to be performed sequentially and uninterruptedly. Although, in many practical cases, this is an acceptable or even standard way of working, this assumption does not allow for the opportunity that maintenance activities of different types are performed separately. Second, the min-cut cut procedure is defined for one maintenance team only. In the current situation at NS this is acceptable; since NS is making a shift to daytime maintenance, it is reasonable to expect that at most one maintenance team is available for daytime maintenance. However, for applications in which the available number of maintenance teams is higher than one, the current min-cut cut generation procedure cannot be used.

Several directions for future research can be recommended.
First, the further development of the \textbf{MSLCP} is considered to be an interesting research area. First, its cut generation process offers opportunities for improvement, and the lessons learned from its development can potentially be used in many other research areas related to scheduling of activities on locations and the capacities of these locations.
Second, for broader applicability, the min-cut cut generation procedure shall be generalized to handle an arbitrary number of teams. This can potentially be achieved by generating additional instants in the \textbf{RAPP}, so that each team has its own dedicated instants. However, additional care must be taken since such an approach may lead to the possibly undesired fact that one maintenance activity can be performed by multiple teams.
Third, an interesting next research topic is the improvement of the computational performance of the \textbf{MLCP}. This improvement is especially relevant in the light of the \textbf{MSLCP} algorithm, since it requires to run the \textbf{MLCP} in each iteration. Looking at the structure of the \textbf{MLCP}, it may potentially be decoupled into multiple smaller sub-problems that are easier to solve, e.g. by decoupling by rolling stock unit, creating sub-problems for each individual rolling stock unit, or by considering a rolling horizon, first optimizing a few days ahead and iteratively adding more days to the optimization.


%% file: main.bbl
\begin{thebibliography}{23}
\providecommand{\natexlab}[1]{#1}
\providecommand{\url}[1]{\texttt{#1}}
\expandafter\ifx\csname urlstyle\endcsname\relax
  \providecommand{\doi}[1]{doi: #1}\else
  \providecommand{\doi}{doi: \begingroup \urlstyle{rm}\Url}\fi

\bibitem[Andr{\'{e}}s et~al.(2015)Andr{\'{e}}s, Cadarso, and
  Mar{\'{i}}n]{Andres2015MaintenanceNetworks}
J.~Andr{\'{e}}s, L.~Cadarso, and A.~Mar{\'{i}}n.
\newblock {Maintenance scheduling in rolling stock circulations in rapid
  transit networks}.
\newblock \emph{Transportation Research Procedia}, 10\penalty0 (July):\penalty0
  524--533, 2015.
\newblock ISSN 23521465.
\newblock \doi{10.1016/j.trpro.2015.09.006}.
\newblock URL \url{http://dx.doi.org/10.1016/j.trpro.2015.09.006}.

\bibitem[Benders(1962)]{Benders1962}
J.~Benders.
\newblock Partitioning procedures for solving mixed-variables programming
  problems.
\newblock \emph{Numerische Mathematik}, 4:\penalty0 238--252, 1962.

\bibitem[Canca and Barrena(2018)]{Canca2018TheSystems}
D.~Canca and E.~Barrena.
\newblock {The integrated rolling stock circulation and depot location problem
  in railway rapid transit systems}.
\newblock \emph{Transportation Research Part E: Logistics and Transportation
  Review}, 109\penalty0 (May 2017):\penalty0 115--138, 2018.
\newblock ISSN 13665545.
\newblock \doi{10.1016/j.tre.2017.10.018}.
\newblock URL \url{https://doi.org/10.1016/j.tre.2017.10.018}.

\bibitem[Clarke et~al.(1997)Clarke, Johnson, Nemhauser, and
  Zhu]{clarke1997aircraft}
L.~Clarke, E.~Johnson, G.~Nemhauser, and Z.~Zhu.
\newblock The aircraft rotation problem.
\newblock \emph{Annals of Operations Research}, 69:\penalty0 33--46, 1997.

\bibitem[Cormen et~al.(2009)Cormen, Leiserson, Rivest, and
  Stein]{cormen2009introduction}
T.~H. Cormen, C.~E. Leiserson, R.~L. Rivest, and C.~Stein.
\newblock \emph{Introduction to algorithms}.
\newblock MIT press, 2009.

\bibitem[Dinmohammadi et~al.(2016)Dinmohammadi, Alkali, and
  Shafiee]{dinmohammadi2016risk}
F.~Dinmohammadi, B.~Alkali, and M.~Shafiee.
\newblock A risk-based model for inspection and maintenance of railway rolling
  stock.
\newblock In \emph{European Safety and Reliability Conference (ESREL)}, pages
  2016--29, 2016.

\bibitem[Feo and Bard(1989)]{Feo1989FlightPlanning}
T.~A. Feo and J.~F. Bard.
\newblock {Flight Scheduling and Maintenance Base Planning}.
\newblock \emph{Management Science}, 35\penalty0 (12):\penalty0 1415--1432,
  1989.
\newblock ISSN 0025-1909.
\newblock \doi{10.1287/mnsc.35.12.1415}.

\bibitem[Ford and Fulkerson(1956)]{fordfulkerson1956}
L.~Ford and D.~Fulkerson.
\newblock {Maximal Flow through a Network}.
\newblock \emph{Canadian Journal of Mathematics}, 8:\penalty0 399--404, 1956.

\bibitem[Gopalan(2014)]{Gopalan2014TheProblem}
R.~Gopalan.
\newblock {The aircraft maintenance base location problem}.
\newblock \emph{European Journal of Operational Research}, 236\penalty0
  (2):\penalty0 634--642, 2014.
\newblock ISSN 03772217.
\newblock \doi{10.1016/j.ejor.2014.01.007}.
\newblock URL \url{http://dx.doi.org/10.1016/j.ejor.2014.01.007}.

\bibitem[Gopalan and Talluri(1998)]{gopalan1998aircraft}
R.~Gopalan and K.~T. Talluri.
\newblock The aircraft maintenance routing problem.
\newblock \emph{Operations Research}, 46\penalty0 (2):\penalty0 260--271, 1998.

\bibitem[Hagberg et~al.(2008)Hagberg, Swart, and S~Chult]{hagberg2008exploring}
A.~Hagberg, P.~Swart, and D.~S~Chult.
\newblock Exploring network structure, dynamics, and function using networkx.
\newblock Technical report, Los Alamos National Lab.(LANL), Los Alamos, NM
  (United States), 2008.

\bibitem[Herr et~al.(2017)Herr, Nicod, Varnier, Zerhouni, and
  Malek~Cherif]{Herr2017F.Scheduling}
N.~Herr, J.-M. Nicod, C.~Varnier, N.~Zerhouni, and F.~Malek~Cherif.
\newblock {F., Joint optimization of train assignment and predictive
  maintenance scheduling}.
\newblock In \emph{7th International Conference on railway operations modelling
  and Analysis (RailLille 2017)}, number April, pages 699--708, 2017.

\bibitem[Hooker(2011)]{hooker2011logic}
J.~Hooker.
\newblock \emph{Logic-based methods for optimization: combining optimization
  and constraint satisfaction}, volume~2.
\newblock John Wiley \& Sons, 2011.

\bibitem[{International Union of Railways}(2018)]{UIC2018paxkm}
{International Union of Railways}.
\newblock Railways, passengers carried, 2018.
\newblock https://data.worldbank.org/indicator/IS.RRS.PASG.KM . Retrieved on
  23-1-2020.

\bibitem[Kravchenko and Werner(2009)]{kravchenko2009minimizing}
S.~A. Kravchenko and F.~Werner.
\newblock Minimizing the number of machines for scheduling jobs with equal
  processing times.
\newblock \emph{European Journal of Operational Research}, 199\penalty0
  (2):\penalty0 595--600, 2009.

\bibitem[Mar{\'{o}}ti and Kroon(2007)]{Maroti2007MaintenanceModel}
G.~Mar{\'{o}}ti and L.~Kroon.
\newblock {Maintenance routing for train units: The interchange model}.
\newblock \emph{Computers and Operations Research}, 34\penalty0 (4):\penalty0
  1121--1140, 2007.
\newblock ISSN 03050548.
\newblock \doi{10.1016/j.cor.2005.05.026}.

\bibitem[Sarac et~al.(2006)Sarac, Batta, and Rump]{Sarac2006ARouting}
A.~Sarac, R.~Batta, and C.~M. Rump.
\newblock {A branch-and-price approach for operational aircraft maintenance
  routing}.
\newblock \emph{European Journal of Operational Research}, 175\penalty0
  (3):\penalty0 1850--1869, 2006.
\newblock ISSN 03772217.
\newblock \doi{10.1016/j.ejor.2004.10.033}.

\bibitem[Taha(2011)]{taha2011operations}
H.~A. Taha.
\newblock \emph{Operations research: an introduction}, volume 790.
\newblock Pearson/Prentice Hall Upper Saddle River, NJ, USA, 2011.

\bibitem[T{\"{o}}nissen and Arts(2018)]{Tonissen2018EconomiesStock}
D.~D. T{\"{o}}nissen and J.~J. Arts.
\newblock {Economies of scale in recoverable robust maintenance location
  routing for rolling stock}.
\newblock \emph{Transportation Research Part B: Methodological}, 117:\penalty0
  360--377, 2018.
\newblock ISSN 01912615.
\newblock \doi{10.1016/j.trb.2018.09.006}.

\bibitem[T{\"{o}}nissen et~al.(2019)T{\"{o}}nissen, Arts, and
  Shen]{Tonissen2019MaintenanceUncertainty}
D.~D. T{\"{o}}nissen, J.~J. Arts, and Z.-J.~M. Shen.
\newblock {Maintenance Location Routing for Rolling Stock Under Line and Fleet
  Planning Uncertainty}.
\newblock \emph{Transportation Science}, 53\penalty0 (5):\penalty0 1252--1270,
  2019.
\newblock ISSN 0041-1655.
\newblock \doi{10.1287/trsc.2018.0866}.

\bibitem[Wagenaar and Kroon(2015)]{Wagenaar2015MaintenanceRailways}
J.~Wagenaar and L.~G. Kroon.
\newblock {Maintenance in Railway Rolling Stock Rescheduling for Passenger
  Railways}.
\newblock \emph{SSRN Electronic Journal}, pages 1--38, 2015.
\newblock \doi{10.2139/ssrn.2566835}.

\bibitem[Zomer et~al.(2020{\natexlab{a}})Zomer, Bešinović, de~Weerdt, and
  Goverde]{GithubCMLCP2020}
J.~Zomer, N.~Bešinović, M.~M. de~Weerdt, and R.~M.~P. Goverde.
\newblock {Demo code for the Maintenance Scheduling and Location Choice Problem
  model}, 2020{\natexlab{a}}.
\newblock GitHub repository. https://github.com/jzomergit/MSLCP/.

\bibitem[Zomer et~al.(2020{\natexlab{b}})Zomer, Bešinović, de~Weerdt, and
  Goverde]{zomer2020MLCP}
J.~Zomer, N.~Bešinović, M.~M. de~Weerdt, and R.~M.~P. Goverde.
\newblock {The Maintenance Location Choice Problem for Railway Rolling Stock},
  2020{\natexlab{b}}.
\newblock arXiv preprint. https://arxiv.org/abs/2012.04565.

\end{thebibliography}
